\documentclass{amsart}

\usepackage[UKenglish]{babel}
\usepackage[UKenglish]{isodate}
\usepackage[T1]{fontenc}

\usepackage{amsmath,amssymb,amsthm}
\usepackage{enumerate}
\usepackage{graphicx}
\usepackage{float}
\usepackage[hidelinks]{hyperref}
\usepackage{color}
\usepackage{todonotes}
\usepackage[style=authoryear-icomp,ibidtracker=false]{biblatex}

\newcommand{\T}{\mathcal{T}}
\newcommand{\M}{\mathcal{M}}
\newcommand{\C}{\mathcal{C}}
\newcommand{\TT}{\mathbb{T}}
\newcommand{\mean}{\operatorname{mean}}
\newcommand{\pr}{\operatorname{pr}}
\newcommand{\BHV}{\mathrm{BHV}}
\newcommand{\CAT}{\operatorname{CAT}}
\newcommand{\rt}{\operatorname{rt}}
\newcommand{\NNI}{\mathrm{NNI}}
\newcommand{\mrca}{\operatorname{mrca}}

\renewcommand{\t}{\mathrm{t}}
\renewcommand{\Im}{\operatorname{Image}}
\renewcommand{\S}{\mathcal{S}}

\newtheorem{theorem}{Theorem}
\newtheorem{proposition}[theorem]{Proposition}
\newtheorem{lemma}[theorem]{Lemma}

\theoremstyle{definition}
\newtheorem{definition}[theorem]{Definition}
\newtheorem{problem}[theorem]{Problem}
\newtheorem{example}[theorem]{Example}

\begin{document}

\title{The space of ultrametric phylogenetic trees}

\author{Alex Gavryushkin}
\address{Department of Computer Science, The University of Auckland, New Zealand}
\email{a.gavruskin@auckland.ac.nz}

\author{Alexei J.\ Drummond}
\email{alexei@cs.auckland.ac.nz}

\begin{abstract}
The reliability of a phylogenetic inference method from genomic sequence data is ensured by its statistical consistency.
Bayesian inference methods produce a sample of phylogenetic trees from the posterior distribution given sequence data.
Hence the question of statistical consistency of such methods is equivalent to the consistency of the summary of the sample.
More generally, statistical consistency is ensured by the tree space used to analyse the sample.

In this paper, we consider two standard parameterisations of phylogenetic time-trees used in evolutionary models: inter-coalescent interval lengths and absolute times of divergence events.
For each of these parameterisations we introduce a natural metric space on ultrametric phylogenetic trees.
We compare the introduced spaces with existing models of tree space and formulate several formal requirements that a metric space on phylogenetic trees must possess in order to be a satisfactory space for statistical analysis, and justify them.
We show that only a few known constructions of the space of phylogenetic trees satisfy these requirements.
However, our results suggest that these basic requirements are not enough to distinguish between the two metric spaces we introduce and that the choice between metric spaces requires additional properties to be considered.
Particularly, that the summary tree minimising the square distance to the trees from the sample might be different for different parameterisations.
This suggests that further fundamental insight is needed into the problem of statistical consistency of phylogenetic inference methods.
\end{abstract}

\maketitle

\section{Introduction}

This paper lies in the broad scope of research on the following two phylogenetic problems, which are also of more general interest, as we demonstrate in this work.
First is the problem of introducing a satisfactory parameterisation of phylogenetic trees for statistical analysis of tree space.
As pointed out by \textcite{Feragen2013-df}, the uniqueness of shortest paths in the space is a desirable property for various types of statistical analysis, while the vast majority of known tree parameterisations do not have this property.
The space of phylogenetic trees encapsulates the structure of a manifold as well as the combinatorially complicated discrete structure of trees \autocite{steel}.
This mix of a continuous and a discrete component is what makes statistical analysis of the space complicated.
The second problem is the problem of summarising a finite set of phylogenetic trees \autocite{heled,hillis2005analysis,huggins}.
This problem arises in different settings of phylogenetic analysis, e.g.\ for computing a statistically consistent summary of a sample from the posterior probability distribution over trees \autocite{beast,beast2}.

An extensive amount of research has been done on the space of phylogenetic trees in the general setting when the phylogenetic distance between taxa is given by arbitrary lengths of the edges of the tree \autocite{bhv,steel}.
As we demonstrate in this paper, this general setting sometimes leads to computationally intractable models when applied to the space of ultrametric trees (a special case of time-trees).
Ultrametric trees are the only satisfactory model for a great body of research in phylogenetics and epidemiology, especially when divergence time dating is the objective, and the taxa are all contemporaneous.
In this case the {\it time-tree} is ultrametric, and is considered separately to the rates of evolution across lineages, which may vary from one branch to the next.

The aim of this paper is to introduce a mathematically satisfactory model of the space of ultrametric phylogenetic trees.
The notion of a `mathematically satisfactory model' will be clarified and made exact later in the paper with an eye towards the two general problems described above.
Our work is inspired by that of \textcite{bhv}, and is similar to it in the sense that we use polyhedral complexes to define a metric space.
The investigation of the tree space from a geometric point of view was initiated by the work \autocite{bhv} with the introduction of a parameterisation that later became known as $\BHV$.
Due to several nice geometric and algorithmic \autocite{owen2011fast} properties, it was recently suggested \autocite{benner2014point} that $\BHV$ is {\em the} space for statistical, and particularly MCMC, analysis of phylogenetic trees.
Our results presented in this paper show how crucial the way a tree is parameterised can be for geometric, algorithmic, and statistical properties of the space.
Particularly, we demonstrate that the summary tree that is suggested in \autocite{benner2014point} will be different for different parameterisations of the tree space.
The question of which parameterisation should be chosen remains open.

Unless otherwise explicit, by a {\em tree} we mean an ultrametric phylogenetic tree, that is, a binary rooted tree with distinguished tips and branch lengths such that the distance from the root is the same to every tip.

We note that although we exclusively consider ultrametric trees in this paper, one of the parameterisations we introduce ($\t$-space) can be generalised to the class of all time-trees as well as to the even more general class of all sampled ancestor trees \autocite{gavryushkina2013recursive,gavryushkina2014bayesian}.
This generalisation is a subject for the future work.

We follow books \autocite{steel} for phylogenetics and \autocite{thomas,bridsonBook} for geometric combinatorics terminology.

\section{Preliminaries}\label{introduction}

It is a standard practice in evolutionary biology to model real biological processes by mathematical abstractions \autocite{steel}.
Particularly, as the goal is often to compare different hypotheses about an evolutionary process modelled by phylogenetic trees, it is natural to work within the space of such trees.
It is also a common practice to introduce different types of measures on the space of trees as a formal way of comparing them.
One of the most general and commonly used ways of measuring the similarity between two trees is given by the notion of a distance, or metric as it is widely known in mathematics.
In order to measure the distance between trees, the tree space has to be {\em parameterised}, that is, some real-valued parameters have to be assigned to trees.

Formally, this scenario can be described as follows.
Let $\T$ be the space of phylogenetic trees on $n$
taxa\footnote{We use $n$ to denote the number of taxa throughout the paper.}.
A {\em parameterisation} of the space $\T$ is an embedding $p\colon \T\to\M$ of the tree space $\T$ to a metric space $\M$, which we call a {\em model metric space}.
By embedding here, we mean a function that maps different trees to different points of the metric space $\M$.
The embedding $p$ plays the role of the assignment of parameters (points of the space $\M$, which could be tuples of real numbers, for example) and the space $\M$ is the parameter space.
The existence of such an embedding makes space $\T$ itself a metric space.
Indeed, the distance between two trees $T$ and $R$ is given by the distance between their images under the embedding $p$, that is, $d_\T(T,R)$ is defined to be $d_\M(p(T),p(R))$.
We say in this case that the metric $d_\T$ is {\em induced} by parameterisation $p$.

As is known \autocite{heled,hillis2005analysis,huggins}, the existence of a parameterisation alone is already a fruitful property of the tree space, as it allows to test hypotheses such as how far are two trees from each other?
How far is an estimate from the true tree?
Given two algorithms, which one produces trees that are closer to the true tree?
Sometimes it is even possible to extract an objective function minimisation that leads to a practical way of summarising posteriors \autocite{heled}.
We will present some of these parameterisations later in this section.

Often phylogenetic analysis requires more subtle properties of the space of trees to be considered, such as what tree is in the middle between two given trees?
What is the path from one tree to
another\footnote{Since we are aimed at a metric space that mirrors the prior or the posterior and preferably both, these two questions are important for us.}?
What is the mean and the variance of a set of (sampled) trees?
The last question is of prominent importance, as this is the very basic question for statistical analysis of data that produces a set of phylogenetic trees.
Furthermore, this question is important in testing whether two probability distributions on tree space are the same, a task common in statistical model selection.
More sophisticated questions include, for example, how standard phylogenetic models such as coalescent and birth-death can be described under a given parameterisation?
Can more efficient proposal mechanisms, such as Hamiltonian Monte Carlo, be employed in Bayesian analysis of phylogenetic data?

A more detailed mathematical analysis is needed in order to approach questions such as these.
In what follows, we summarise several basic properties of parameterisations, which we suggest are desirable to advance research on the problems mentioned.

It is often the case that the metric space $\M$, that is used to parameterise tree space $\T$, is greatly different from the metric space $\T$ with the induced metric $d_\T$.
The key reason for this is the nature of the parameterisation $p$.
As we will see later in the paper, some parameterisations $p$ induce metrics that share almost no geometric properties in common with the original metric space $\M$ that was used in the parameterisation $p$.
Particularly, those parameterisations are far from being {\em bijective}, that is, being able to recover a tree given an arbitrary point from the space $\M$.
The lack of this property can lead to situations where, for example, there are infinitely many trees all of which minimise the total square distance to a given set of trees \autocite{heled}.

Although the parameterisations we introduce in Sections~\ref{tauSpace} and~\ref{tSpace} of this paper are bijective, the requirement of being bijective is somewhat strong in the sense that many desirable properties can be achieved without the parameterisation being bijective.
We continue with introducing formal requirements that allow to carry the analysis of the space $\M$ over to the space of trees $\T$.

For the statistical analysis of a space, one needs to define probability distributions over the space, e.g.\ for Bayesian analysis the first step is to define a prior distribution.
A continuous probability distribution defined on the metric space $\M$ has to remain the
same\footnote{In the sense that all statistics, e.g.\ $k$-th moments, are preserved.}
continuous distribution when pulled back to the space of trees $\T$ under the parameterisation $p$.
In order to achieve this, one has to be able to continuously move from one tree to another by a path that stays within the tree space.
In other words, any two trees have to be connected by a path.

Formally, a metric space $X$ is called {\em path-connected} if for each pair of points $x,y$ in the space, there exists a continuous map $\gamma$ (with respect to the standard topologies generated by balls) from the unit real segment $[0,1]$ to the space $X$ such that $\gamma(0)=x$ and $\gamma(1)=y$.

Thus, the first property a satisfactory parameterisation of the tree space must satisfy is:
\begin{equation}\label{pathconnected}
\Im(p)\mbox{ is path-connected in $\M$.}\tag{P1}
\end{equation}

Our next property ensures that (shortest) paths in model metric space $\M$ remain (shortest) paths when pulled back to tree space $\T$. A subspace $X$ of a space $Y$ is called {\em convex} if for every pair of points $x,y\in X$, every shortest path $\gamma$ between $x$ and $y$, and every real number $s\in[0,1]$, it follows that $\gamma(s)\in X$.
\begin{equation}\label{convex}
\Im(p)\mbox{ is convex in $\M$.}\tag{P2}
\end{equation}

The next requirement is necessary to specify a probability distribution over trees by defining a probability distribution over the model metric space.
For this method to work, the $\Im(p)$ has to be a non-trivial part of $\M$:
\begin{equation}\label{dimension}
\Im(p)\mbox{ has the same dimension as $\M$.}\tag{P3}
\end{equation}

Requirements~\ref{pathconnected}--\ref{dimension} guarantee that desirable geometric properties of space $\M$ will be inherited by the induced metric space on trees $\T$, but none of the requirements causes those properties to exist.
They have to be postulated.
Hence, we now go on to the properties of the space $\M$.
It is important to note that the following properties only make sense if the requirements~\ref{pathconnected}--\ref{dimension} are fulfilled.

Our next requirement has to do with the uniqueness of shortest paths, which is a necessary property for statistical analysis \autocite{Feragen2013-df}.
The uniqueness of shortest paths implies the uniqueness of several types of means, the soundness of the notion of a variance, and the existence and uniqueness of summary trees obtained by minimising an objective function of square distance.

We say that a metric space {\em possesses unique geodesics} if there exists a unique shortest path between every two points in the space.
This shortest path is called a
{\em geodesic}\footnote{It is worth noting here that our notion of geodesic is somewhat different from the one that is commonly used in differential geometry.
We call a path geodesic only if the path is {\em globally} shortest.
For example, the great circle of a sphere with a small interval removed is {\em not} a geodesic in our sense.}.
\begin{equation}\label{geodesic}
\mbox{Metric space $\M$ possesses unique geodesics.}\tag{P4}
\end{equation}

This requirement can in practice be relaxed to hold almost surely.
Intuitively this means that with probability one the shortest path is unique between two points drawn at random.
Formally, we assume that the metric space $\M$ is equipped with a probability measure $\mu$ and say that a property $P(\cdot)$ is satisfied {\em almost surely} if $\mu^* \{ x \in \M \mid P(x) \} = 1$, where $\mu^*$ is the product measure if property $P$ is defined on tuples.
The relaxed requirement is:
\begin{equation}\label{geodesic_almost}
\mbox{Metric space $\M$ possesses unique geodesics almost surely.}\tag{P$4^\prime$}
\end{equation}

A sphere with the standard spherical distance and uniform measure gives an example distinguishing properties~\ref{geodesic} and~\ref{geodesic_almost}.

Since geodesics can be incomputable for some metric
spaces\footnote{It is not hard to see that the halting problem for Turing machines can be reduced to the problem of computing shortest paths in graphs.
More precisely, there exists a computable graph $G$ such that any algorithm that computes shortest paths between vertices in $G$, solves the halting problem.},
our next property of model space $\M$ is:
\begin{equation}\label{computable}
\mbox{Geodesics in metric space $\M$ are computable.}\tag{P5}
\end{equation}

A natural strengthening of Property~\ref{computable} that is necessary to make the parameterisation potentially useful in practice is:
\begin{equation}\label{efficient}
\mbox{Geodesics in metric space $\M$ are efficiently computable.}\tag{P$5^\prime$}
\end{equation}

The computational complexity of geodesics is fundamental for applications, as the algorithms for computing various characteristics of a data set such as the mean, variance, diversity, confidence regions, and so on rely on computing geodesics as a subroutine \autocite{bacak2012computing,owen2011fast}.

Our work is motivated by the lack of parameterisations in the literature that enjoy all properties~\ref{pathconnected}--\ref{computable}.
Indeed, all known summary tree estimators operate in spaces larger than the space of ultrametric rooted binary trees, hence breaking requirement~\ref{dimension}.
For instance, \textcite{heled} and \textcite{huggins}
%this is paper [8] from Joseph's paper
use the so-called Rooted Branch Score (RBS) metric space for producing a summary tree given a sample of trees from the posterior distribution.
The idea of the RBS space is to encode a tree on $n$ taxa by a $(2^n-1)$-dimensional real vector, find an optimum in the $(2^n-1)$-dimensional Euclidean space, and find the nearest point in the Euclidean space that can be pulled back to the tree space.
Although this approach proved to be fruitful in several applied scenarios \autocite{heled}, it lacks properties~\ref{pathconnected}--\ref{dimension}.
Moreover, a tree that minimises the RBS distance to a (finite) set of trees is not unique---indeed there could be infinitely many such trees.
This optimisation problem is computationally intractable even for moderate values of $n$.
In implementations of this method, the inefficiency is overcome by restricting the search only to tree topologies that are present in the posterior sample, that is, in the given set of trees.
Furthermore, the tree topologies and the branch lengths have to be summarised separately in order to make the method computationally tractable \autocite{heled}.

Other metrics used by \textcite{huggins} employ projections to smaller dimension spaces to overcome the absence of properties~\ref{pathconnected}--\ref{dimension}.
Those metrics share the same pathologies as RBS.
Moreover, the use of projections for estimating means can lead to unbounded errors as witnessed by the following proposition that claims that the projection of the mean can be as far from the mean of the projections as possible.

\begin{proposition}
Let $N$ be a (arbitrarily large) real number, $E$ a Euclidean space of dimension $k>1$, and $x_1,\ldots,x_s$ a set of points in $E$.
Then there exists a subspace $D$ of $E$ such that
\[
d_E(\pr_D(\mean_E(x_1,\ldots,x_s)),\mean_D(\pr_D(x_1),\ldots,\pr_D(x_s)))
\geq N,
\]
where $d_E$ is the Euclidean distance, $\pr_D(x)$ is the projection of the point $x\in E$ onto $D$, and $\mean_X(x_1,\ldots,x_s)$ is the Fr\'echet mean of $x_1,\ldots,x_s$ in the space $X$.
\end{proposition}

\proof
We prove the proposition for $k=s=2$.
An arbitrary case is analogous.
Let $\ell$ be the line through $x_1$ and $x_2$ in $E$ and $\ell_0$ be a line parallel to $\ell$ at a distance $M$ from $\ell$.
Consider a parabola $D$ which has its vertex on the line $\ell_0$ and crosses the line $\ell$ at some points $a$ and $b$ both of which are between $x_1$ and $x_2$.
It is not hard to see that for large enough $M$, we get $d_E(\pr_D(\mean_E(x_1,x_2)),\mean_D(\pr_D(x_1),\pr_D(x_2)))\geq N$.
\endproof

It might appear that the construction used in the proof is artificial, but this is actually very similar to what is happening in such parameterisations as RBS and dissimilarity map distance \autocite{huggins}, where the conditions on the set of points that correspond to trees are non-trivial \autocite{cardona2010nodal}.
The dissimilarity map distance \autocite{huggins} between two trees is defined as the distance between the distance matrices of the trees, in the space of square matrices.
That is, the parameterisation $p$ maps a tree to its distance matrix, and the model metric space $\M$ is the space of $n\times n$ matrices with the pointwise distance.
This space is geometrically similar to RBS in the way that none of properties~\ref{pathconnected}--\ref{dimension} are satisfied.
\textcite{cardona2010nodal} characterised $\Im(p)$ for the case when the trees are not necessarily ultrametric.
This characterisation fulfils the requirements~\ref{pathconnected}--\ref{dimension}.
An attempt to carry this characterisation over to the space of ultrametric trees has the same complication as BHV space, which we discuss below.

The most geometrically attractive parameterisation of the (non-ultrametric) tree space is the $\BHV$ space \autocite{bhv}.
This is the only parameterisation we are aware of that fulfils all the properties~\ref{pathconnected}--\ref{efficient} \autocite{bhv,owen2011fast}.
This parameterisation employs a $(2n-2)$-dimensional cubical complex with unique geodesics as the model metric space $\M$, then a bijective correspondence between the space of all phylogenetic trees and the complex $\M$ is established.
Trees of a fixed topology are parameterised by a $(2n-2)$-dimensional vector given by the lengths of the branches, and correspond to a cube.
The adjacent cubes of the complex correspond to NNI-adjacent trees.
Although it took ten years to establish property~\ref{efficient} for the parameterisation, the polynomial algorithm designed by \textcite{owen2011fast} appears to be quite practical.

As we demonstrate in the next section, it is somewhat involved to apply the BHV model, as well as other BHV-like models \autocite{Feragen2013-df, Miller2015-rk}, to the space of ultrametric trees.

\section{Preliminary attempt}\label{preAttemptSection}

A possible (naive) approach could be to simply restrict the $\BHV$ space to the set of ultrametric trees.
Unfortunately, this simple adaptation lacks all properties~\ref{pathconnected}--\ref{dimension}, so the algorithms developed by \textcite{owen2011fast} become inapplicable.

Another (less naive) approach is to parameterise a tree by the lengths of all internal edges and the shortest external edge.
In this case, the lengths of the rest of external edges are computed so that the resulting tree is ultrametric.
This `less naive adaptation' of $\BHV$ space is similar to the `bounded $\BHV$' adaptation, which we consider later in this section.

A fundamental characteristic of all $\BHV$-like spaces is that the subspaces corresponding to different ranked tree topologies have different volumes.
This property results in complications for introducing a (prior) probability distribution over the space.

In the rest of this section, we model the space of trees by a set of bounded polyhedral complexes indexed by the set of positive reals.
We assume here that the reader is familiar with $\BHV$ space \autocite{bhv}.
Otherwise, the rest of this section (excluding the next paragraph) can be skipped, as the following sections of the paper are self-containing.

Since the complexity of presentations is not the matter of this paper, we shall make no distinction between the tree space $\T$ and the model metric space $\M$ used in the parameterisation $p$ of $\T$, in the case when $p$ is a bijection.
For instance, when $\M$ has unique geodesics and $p$ is a bijection, we shall simply say that $\T$ has unique geodesics (under this parameterisation).
A parameterisation $p$ is called {\em strict} if $p$ is a bijection.

Consider the space $\BHV^\circ$, which is the $\BHV$ space where external branches are ignored, that is, the projection of $\BHV$ to the coordinates corresponding to internal branches.
We restrict each orthant of space $\BHV^\circ$ to the set $\{T\mid T$ has height at most $H\}$, where $H$ is a fixed real number, and denote thus obtained space by $\BHV^\circ\upharpoonright H$.
Space $\BHV^\circ\upharpoonright H$ can be seen as the space of trees of height $H$ because every tree from $\BHV^\circ\upharpoonright H$ can be extended in a unique way to a tree of height $H$ by attaching the external edges of appropriate lengths to the places where they were in the original $\BHV$ space.
Thus, the polyhedral complex $\BHV^\circ\upharpoonright H$ is a strict parameterisation of the space of ultrametric trees of height $H$.
By varying $H$ over the set of positive reals, we get a strict parameterisation of the tree space as a set of bounded polyhedral complexes indexed by positive reals.
We call this space {\em bounded $\BHV$ space.}

Although the space $\BHV^\circ\upharpoonright H$ is not a cubical complex, it is geometrically and algorithmically similar to the $\BHV$ space.
Indeed, since in a neighbourhood of the origin the space $\BHV^\circ\upharpoonright H$ is a cubical complex, it possesses efficiently computable unique geodesics in the same way as $\BHV$ does.
This can be seen by noticing the following.
Suppose $\C$ is a cubical complex with unique geodesics such that each cube is given by inequalities $x_i\leq K$.
Let $\S$ be a polyhedral complex obtained from $\C$ by replacing the inequalities $x_i\leq K$ by $\Sigma_ix_i\leq K$.
Then $\S$ has unique geodesics.
Furthermore, if geodesics in $\C$ are
efficiently\footnote{By `efficiently computable' here and in the rest of the paper we mean computable in (low degree) polynomial (in the number $n$ of taxa) time.}
computable then so are geodesics in $\S$.
Both of the statements are not hard to prove, but this goes beyond the scope of this paper.

The first and most obvious complication of this parameterisation is the lack of independence between coordinates.
The last coordinate, the height of the tree, cannot be smaller than the sum of coordinates corresponding to the internal edges.
This results in non-trivial boundary conditions that has to be taken into account in the study of the geometry of the space, and more problems with implementing algorithms.
Another feature of this space is that a change of the length of only one internal branch causes a change of the length of all external edges.
Hence, if the edge length is interpreted as time, which is the case for many phylogenetic applications, then a change of an older divergence time impacts the times of most recent divergence events for {\em each} taxon.

More fundamental issues with this parameterisation are the following.
If (some of) the branch lengths are given by confidence intervals then computing the confidence region in the space becomes a non-trivial exercise.
We already mentioned above that the non-uniform distribution of the volume among different ranked tree topologies in the space makes it difficult to introduce (prior) probability distributions used in (Bayesian) inference of time-trees.

To overcome these and similar issues is the goal of this paper.

\section{$\tau$-space}\label{tauSpace}

In this section, we model the space of ultrametric trees by a cubical complex, which we call $\tau$-space, with efficiently computable unique geodesics and establish several geometric and algorithmic properties of the space.

\subsection{Construction of space}

\begin{figure}
\centering
\includegraphics[width=0.56\textwidth]{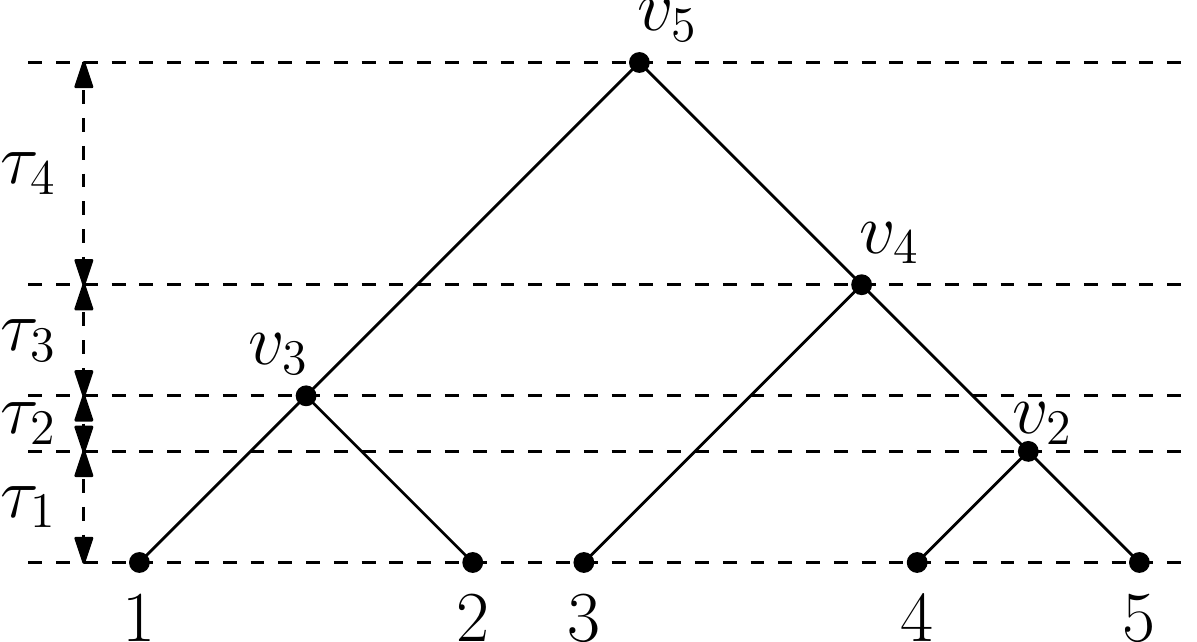}
\caption{Parameterisation of tree from $5$-dimensional $\tau$-space $\T_5$.}
\label{T5.pdf}
\end{figure}

\begin{figure}
\centering
\includegraphics[width=\textwidth]{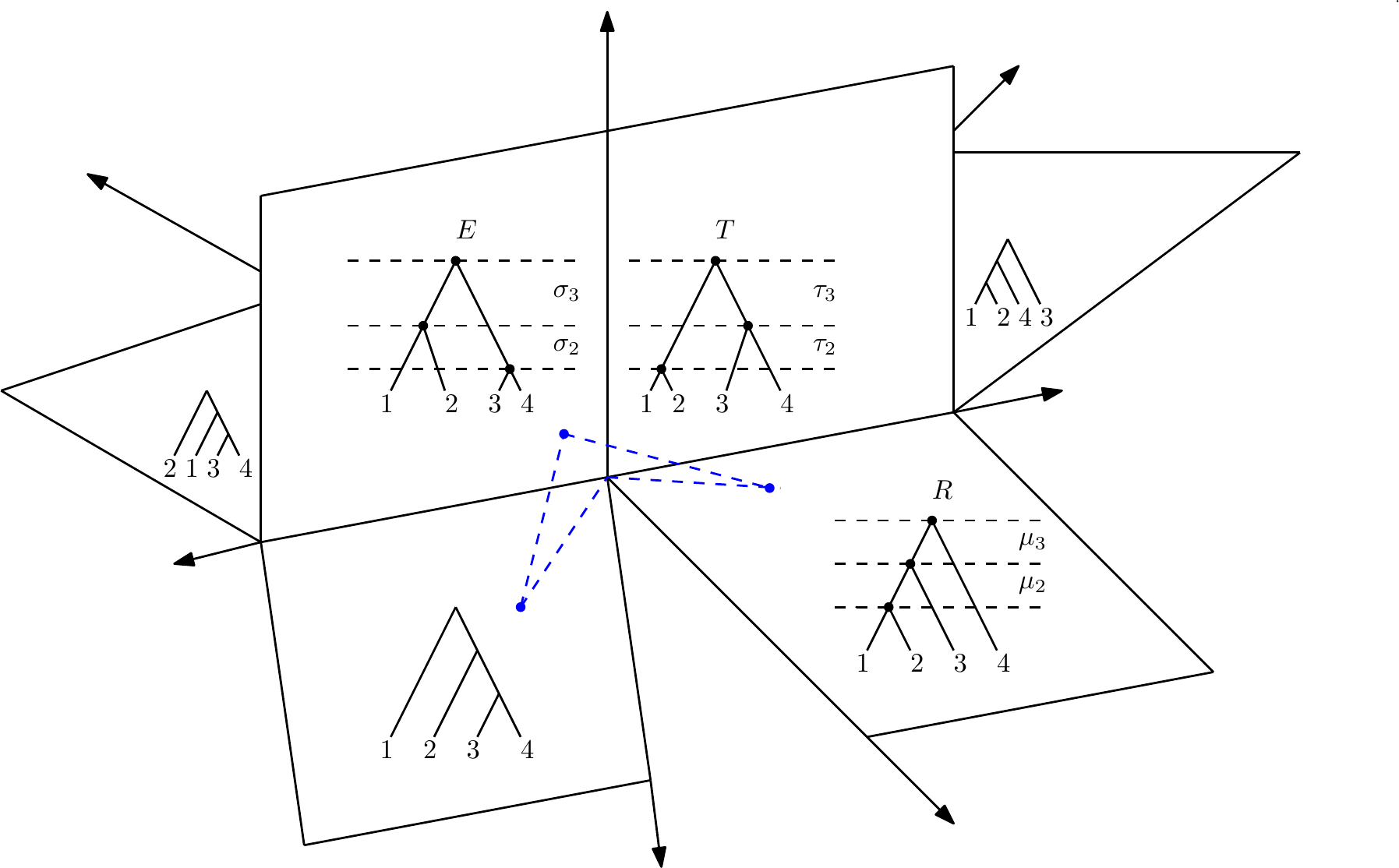}
\caption{Three-dimensional projection of one third of $4$-dimensional $\tau$-space $\T_4$.
Each orthant is projected onto the subspace with the first coordinate $\tau_1$ fixed.
Although the projected space cannot be embedded into $3$-dimensional Euclidean space, it can be visualised by imagining the other two thirds of the space.
Triangles are thin in the space due to the Cartan--Alexandrov--Toponogov axiom (for $k=0$), see Definition~\ref{definitionCAT}.
}
\label{tauSpace.pdf}
\end{figure}

We begin with a formal construction of the space illustrated in Figures~\ref{T5.pdf} and~\ref{tauSpace.pdf}.
We will be using terms `height of a node' and `time of a node' interchangeably to refer to the distance from a taxon to the node.

Let $T$ be an ultrametric tree on $n$ taxa with times assigned to its nodes.
Assuming that the times of all internal nodes are pairwise distinct, we denote the set of such trees by $\T_n$.
We parameterise tree $T$ by a pair that consists of the ranked topology of the tree and the differences between the times of the tree's consecutive nodes.
We proceed by defining this parameterisation in detail.
Let us order the internal nodes of $T$ according to their times: $v_2,\ldots,v_{n}$.
Note that the node $v_{n}$ must be the root in this case.
Denote the difference between the time of node $v_{i+1}$ and the time of node $v_{i}$ by $\tau_{i}$ for all $i\in\{2,\ldots,n-1\}$.
We call $\tau_i$ the {\em coordinate of node} $v_i$.
Since the tree is ultrametric, the differences between the time of $v_2$ and the times of external nodes are all the same.
Denote this difference by $\tau_1$.
The coordinates of tree $T$ are given by the $n$-tuple $(\rt(T),\bar\tau)$, where $\rt(T)$ is the ranked topology of tree $T$ and $\bar\tau$ is the tuple $(\tau_1,\ldots,\tau_{n-1})$ from $\mathbb R_0^{n-1}$ that consists of the coordinates of the nodes of $T$.
By $\mathbb R_0^{n-1}$ we denote the $(n-1)$-dimensional non-negative orthant $\{(r_1,\ldots,r_{n-1})\mid r_i\in\mathbb R~\&~r_i\geq0\}$, where $\mathbb R$ is the set of reals.
Figure~\ref{T5.pdf} depicts an example of $\tau$-parameterisation of a tree from $\T_5$.

Consider now the set of all ranked topologies on $n$ taxa such that all internal nodes have different ranks.
We recall that there are $\frac{(n-1)!\cdot n!}{2^{n-1}}$ many such topologies \autocite{steel}, and we denote this number by $m$ throughout the paper.

Thus, we have constructed a disjoint union of $m$ $(n-1)$-dimensional polyhedra $S=\{(\rt(T),\bar\tau)\mid T\in\T_n,\bar\tau\in\mathbb R_0^{n-1}\}$.
Specifically, the polyhedra are orthants indexed by tree topologies.
It is clear that the set $\T_n$ is in a bijective correspondence with the interior of $S$.
It is also obvious how to establish a bijection between the faces of the polyhedra in $S$ and the set of ranked (multifurcating) tree topologies on $n$ taxa which have at least two internal nodes of the same rank.
Indeed, if we consider such a tree, the coordinates $\tau_i$ that are between two nodes of the same rank have to be $0$, and the faces of the polyhedra in $S$ are precisely the tuples $(\rt(T),\bar\tau)$ where some of the coordinates $\tau_i$ are $0$.

We now want to create a polyhedral complex in the obvious way, that is, by gluing the faces that correspond to same ranked (not necessarily completely resolved) tree topologies together.
We proceed formally as follows.
Let us define an equivalence relation $\sim$ on the set of faces of polyhedra in $S$.
We say that two faces $F$ and $G$ are equivalent, written $F\sim G$, if they correspond to the same ranked tree topology.
Now, consider the set $S'$ that consists of the union of the set $S$ and the set of all faces of elements from $S$.
The polyhedral complex is then the quotient set $S'/\sim$.

Since trees are in a bijective correspondence with this complex, the parameterisation is strict and from now on we shall identify the space of trees $\T_n$ with this polyhedral complex, slightly abusing the
notation\footnote{We note that we abuse the notation here not only because we make no distinction between the space of trees and the polyhedral complex, but also because the multifurcating trees are present in the complex and absent from the tree space $\T_n$ we initially considered.}.

We shall assume that all dimensions of the orthants $\tau_i$ are bounded from above by a (large enough) constant.
This boundary makes the polyhedral complex $\T_n$ a cubical complex, which is a standard and well-studied object of geometric combinatorics \autocite{bridsonBook,thomas}.
This restriction is inessential for this paper as all the results remain true in the unbounded case.

An interesting example where asymptotic properties of tree space at infinity are treated differently is so-called `orange geometry' \autocite{Kim2000-ba}, where polyhedra are `glued' together at infinity.
That geometric framework encodes not only phylogenies but also data sets, evolution models, and
estimation methods.

From now on tree space $\T_n$ is a cubical complex.
We assume Euclidean distance within cubes.
This assumption implies that space $\T_n$ is a metric space with geodesics.
Indeed, since $\T_n$ is a finite connected cubical complex, every path is $\T_n$ can be partitioned into a finite number of subpaths each of which is a Euclidean path.
The length of a path in $\T_n$ is then the sum of the lengths of those subpaths.
The distance between two points in $\T_n$ is given by the length of a shortest path between the points.
Every shortest path between two points is called a {\em geodesic}.
We prove later in this paper that geodesics are (globally) unique in $\T_n$.
Hence, $\T_n$ is a metric space with unique geodesics.
We call this space {\em $\tau$-space}.

\subsection{Geometric properties}\label{tauGeometry}

Let us start with consider geometric properties of the space and compare them with those of $\BHV$ space.
It is convenient to think of $\tau$-space as a set of points (trees) that freely move within cubes without leaving them as long as all the coordinates $\tau_i$ are strictly positive.
This movement results in changes in branch lengths (waiting time between divergence events) but the ranked tree topology remains the same.
When one of the $\tau_i$ becomes $0$, the point (tree) is on the boundary of the cube of one smaller dimension.
The point now can either move along the boundary by varying the other $\tau_i$, or it can leave the boundary by increasing the $\tau_i$ that became $0$.
The boundary corresponds to a
facet\footnote{By a {\em facet} of a polyhedron here and throughout the paper, we mean a face the dimension of which is one smaller than the dimension of the polyhedron, that is, a face of codimension one.}
$F$ and there could be several cubes that share this facet $F$.
It is not hard to understand that the possible numbers of cubes which share a common facet are one, two, and three.
Indeed, setting $\tau_1 = 0$ gives an example of a cube and its facet that is not a face of any other cube.
If a facet does not correspond to a multifurcation (see the facet between cubes corresponding to trees $T$ and $E$ in Figure~\ref{tauSpace.pdf}), there are precisely two cubes that share the facet.
If it does (as the other facet of the cube corresponding to tree $T$ in Figure~\ref{tauSpace.pdf}), then the number is three.

\begin{figure}[H]
\centering
\includegraphics[width=.85\textwidth]{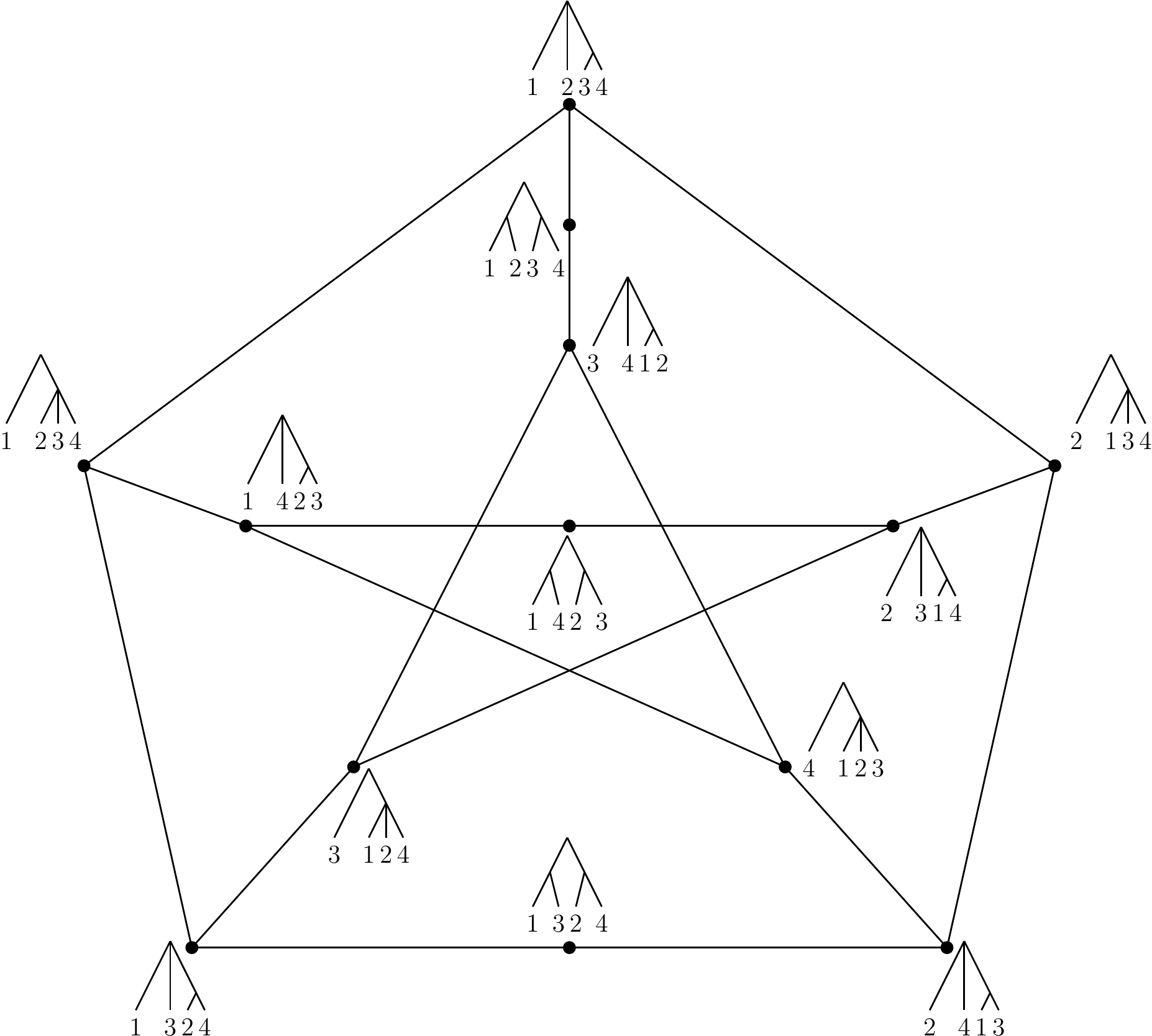}
\caption{Link of origin of $\T_4$.}
\label{tauLink.pdf}
\end{figure}

At first glance it might seem that the $\BHV$ and $\tau$-space are very
similar\footnote{This subsection can be skipped by those who are not familiar with $\BHV$ space, as the rest of the paper does not depend on this subsection.}.
The graph in Figure~\ref{tauLink.pdf} depicts the link of the origin of space $\T_4$.
The link is similar to that of $\BHV$ space on four taxa indeed, but it already suggests several differences that we would like to investigate.

In this subsection, we establish several geometric properties of the two spaces to better understand the differences and similarities between them in order to answer the question of whether the algorithms developed for $\BHV$ space \autocite{owen2011fast} are applicable in $\tau$-space.
The first property we want to point out is the following.
For every tree topology, the dimensions of corresponding orthants in $\BHV_n$ and $\T_n$ are different.
This is because external edges add $n$ to the dimension of the $\BHV$-orthant and add $1$ to the dimension of the $\tau$-orthant.
One might suggest that the spaces $\BHV_n^\circ$ and $\T_n^\circ$, which are the corresponding spaces where all external edges are omitted, are geometrically similar.
They are indeed, they share a number of geometric properties.
However, some dissimilarities between them become clear if one attempts to uniformly map one distance to the other.
If such a mapping existed, all geometric and algorithmic results for $\BHV$ could be directly applied to $\tau$-space.
We formalise this assertion in the following two propositions.

\begin{proposition}
Spaces $\BHV_n^\circ$ and $\T_n^\circ$ are not isometric.
\end{proposition}

\proof
This follows from the fact that isometries preserve angles.
Indeed, let us fix a non-caterpillar tree topology and consider the corresponding orthants in $\BHV_n^\circ$ and $\tau_n^\circ$.
We may notice that there are several orthants corresponding to the tree topology in $\tau_n^\circ$ and only one orthant in $\BHV_n^\circ$ space.
One can use an appropriate number of hyperplanes to partition the $\BHV_n^\circ$-orthant in a way that every member of the partition corresponds to the trees in precisely one $\tau_n^\circ$-orthant.
Clearly, no embedding between these subspaces preserves angles.
\endproof

\begin{figure}[H]
\centering
\includegraphics[width=0.67\textwidth]{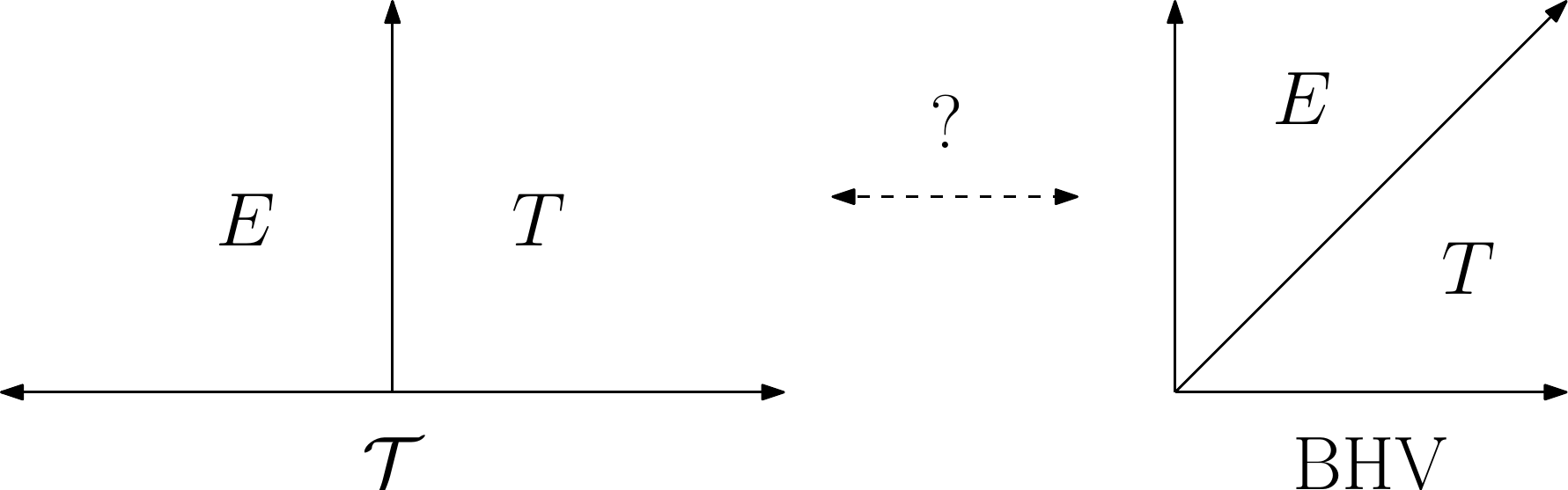}
\caption{$\BHV$ and $\tau$-space are not isometric.}
\label{noIsometry.pdf}
\end{figure}

The proof above can intuitively be understood by trying to establish an isometry between the orthants that correspond to the trees $T$ and $E$ in Figure~\ref{tauSpace.pdf}.
The corresponding $\tau$- and $\BHV$-subspaces can be drawn as Figure~\ref{noIsometry.pdf} (note that the objects depicted are flat).

Another seemingly plausible hypothesis could be that the $\BHV$-distance majorates the $\tau$-distance, that is, $d_\tau(T,R) \leq d_{\BHV}(T,R)$ for all trees $T$ and $R$.
Although it is obvious that the $\BHV$- and the $\tau$-coordinates are easily computable from each
other\footnote{It is important to note that this claim is about {\em coordinates} rather than {\em distances}.
Computing the $\tau$-distance given the $\BHV$-distance, or vice versa, can be somewhat involved in general.}, the following proposition is true.

\begin{proposition}
None of the $\BHV$- and $\tau$-metrics majorates the other.
\end{proposition}

It is important to note that since the dimensions of $\BHV_n$ and $\tau_n$ are different, we are ignoring the external branches here and considering $\BHV_n^\circ$ and $\tau_n^\circ$.

\proof
Consider trees $T$, $R$, and $E$ depicted in Figure~\ref{tauSpace.pdf}.
We finish the proof of the lemma by setting

(1) All sigmas, mus, and taus to 1.
In this case:\\
$d_\tau(T,R) = 2 < \sqrt5 = d_{\BHV}(T,R)$\\
$d_\tau(T,E) = 2 > \sqrt2 = d_{\BHV}(T,E)$

One might hypothesise that an inequality of the second type can only be obtained in the quadrants that present in $\tau$-space but not in $\BHV$ space.
Although it is not necessary for the proof, we demonstrate that this hypothesis can be refuted by setting

(2) $\tau_2=1$, $\tau_3=2$, $\mu_2=3$, $\mu_3=4$.
In this case:\\
$d_\tau(T,R) = \sqrt{40} > 6 = d_{\BHV}(T,R)$.
\endproof

\subsection{Uniqueness and efficiency of geodesics}

The geometric property of main interest to us in this paper is the (global) uniqueness of geodesics, as this property is crucial for such geometric characteristics as the Fr\'echet mean, standard deviation, convex hulls, etc.

We recall that a metric space $X$ is called {\em geodesic} if every pair of points from $X$ is connected by a shortest path.
A geodesic metric space is said to {\em have unique geodesics} if the geodesic is unique between every two points from $X$.

\begin{definition}\label{definitionCAT}
A geodesic metric space $X$ is said to {\em satisfy Cartan--Alexandrov--Toponogov axiom,} or {\em be $\CAT(0)$,} if the following property holds.

\begin{itemize}
\item[]
For all triples $x_1,x_2,x_3\in X$ and all points $y$ on a geodesic from $x_1$ to $x_2$, the inequality $d_X(x_3,y)\leq d_E(x_3',y')$ holds, where $x_1',x_2',x_3'$ are three points on the Euclidean plane such that $d_E(x_i',x_j')=d_X(x_i,x_j)$ for all $1\leq i<j\leq3$, $y'$ is the point on the segment $[x_1',x_2']$ such that $d_E(x_1',y')=d_X(x_1,y)$, and $d_E$ is the Euclidean distance.
\end{itemize}
\end{definition}

In other words, a metric space $X$ is $\CAT(0)$ if no triangle $\Delta$ in $X$ is thicker than a Euclidean triangle $\Delta_E$ of the same size as $\Delta$.

It follows from the definition of a $\CAT(0)$ metric space that the space has unique geodesics.
Indeed, let $X$ be a $\CAT(0)$ space and $a,b$ two points from $X$.
Consider a point $x$ on a geodesic $\gamma$ from $a$ to $b$ and consider a degenerate Euclidean triangle $a',x',b'$ where $x'$ lies on the segment $[a',b']$ at the same distance from $a'$ as $x$ is from $a$ in $X$.
The axiom $\CAT(0)$ implies then that $d_X(x,y)\leq d_E(x',x')$, where $y$ is a point on any geodesic from $a$ to $b$ at the same distance from $a$ as $x$.
Since $d_E(x',x')=0$, $d_X(x,y)=0$ and every geodesic from $a$ to $b$ coincides with $\gamma$ because we have chosen $x$ arbitrarily.

We derive the fact that $\tau$-space has unique geodesics from the following theorem.
Note that a cubical complex is said to have the {\em intrinsic Euclidean metric} if the complex is metrised in the same way as we metrised $\tau$-space.

\begin{theorem}[\cite{gromovOriginal}]
A cubical complex $\C$ with the intrinsic Euclidean metric is $\CAT(0)$ if and only if $\C$ is connected, simply connected, and for all natural numbers $k$, if three $(k+2)$-cubes of $\C$ share a common $k$-cube and pairwise share common
different
%TONOTE: The word `different' does not appear in the original statement, but without it the theorem has a simple counterexample.
$(k+1)$-cubes, then they are contained in a $(k+ 3)$-cube of $C$.
\end{theorem}

Clearly $\tau$-space is a cubical complex which is connected and simply connected (see also Lemma~\ref{homeomorphic} below).
For the last requirement of the theorem we note that the $(k+2)$-cubes cannot be of the highest possible dimension, otherwise the $(k+1)$-cubes would result in a cycle of length $3$ in the link of the
origin\footnote{The {\em link of a vertex} $v$ of a polyhedral complex is defined as a graph with nodes being the facets that contain $v$, where two nodes $x,y$ are adjacent in the graph if there exists a polyhedron in the complex with facets $x$ and $y$.}
of $\tau$-space.
The impossibility of such a cycle follows from the fact that the $\NNI$ graph does not have $3$-cycles and also can be shown in the same way as we consider longer cycles in the proof of Theorem~\ref{tSpaceCAT0}.
Hence, we can assume that there exists a
unique\footnote{The uniqueness is assumed for the sake of clarity.
The proof can straightforwardly be modified to the case when there are several non-resolved $\tau$-coordinates.}
$\tau_i$ such that the first $(k+2)$-cube $C_1$ has $\tau_i = 0$ and has the rest of $\tau$-coordinates strictly positive.
Similarly, the second $(k+2)$-cube $C_2$ has a unique $\tau_j = 0$ and $C_3$ has $\tau_r = 0$.

\begin{itemize}
\item[Case 0.] $i = j = r$.
This case results in a cycle of length $3$ in the link of the origin of $\tau$-space, which is impossible.
\item[Case 1.] $i \ne j = r$.
Since $C_2$ and $C_3$ share a $(k+1)$-cube, they both must have a coordinate $\tau_s$ such that $\tau_s > 0$ and if we set $\tau_s = 0$ in both $C_2$ and $C_3$ then the resulting cubes coincide.
We note that $s \ne i$ because the $(k+1)$-cubes must be pairwise different.
Then the only way for $C_1$ and $C_2$ to share a $(k+1)$-cube is via setting both $\tau_i$ and $\tau_j$ to zero (this is because $i \ne j$).
Hence $\tau_s$ is resolved in the same way in $C_1$ and $C_2$.
By the same reason, $\tau_s$ is resolved in the same way in $C_1$ and $C_3$.
This implies that $C_2$ and $C_3$ coincide, so this case is impossible.
\item[Case 2.] All $i$, $j$, and $r$ are pairwise distinct.
Since $C_1$ and $C_2$ share a $(k+1)$-cube, $\tau_r$ is resolved in the same way in both of these cubes.
By a similar reason, $\tau_i$ is resolved in the same way in $C_2$ and $C_3$ and $\tau_j$ in $C_1$ and $C_3$.
In this case we construct a $(k+3)$-cube that contains all $C_1$, $C_2$, and $C_3$ by taking cube $C_1$ and resolving $\tau_i$ in the way it is resolved in $C_2$ and $C_3$.
\end{itemize}

Thus, we have established the following result.

\begin{theorem}\label{uniqueGeodesics}
$\tau$-space has unique geodesics.
\end{theorem}

This property is fundamental for summarising sets of trees, because the uniqueness of geodesics implies that several geometric centres are unique.
For example, such objects as Fr\'echet mean \autocite{Karcher1977-bu}, standard deviation, convex hull, and many other, are well-defined.
Furthermore, since the notions of a mean and a variance are well-defined, fundamental theorems of probability theory, such as the Central Limit Theorem, can be studied in tree space \autocite{Barden2013-ov, Miller2015-rk, Nye2011-aa, Nye2015-re}.

The next question we would like to study is the question of effectiveness.
Can different types of means, centroids, and hulls be efficiently computed?
The answer to this question depends on whether or not geodesics are efficiently computable.
Indeed, in most cases the non-existence of a polynomial algorithm for computing geodesics implies the non-existence of such algorithms for computing various versions of means.

A careful consideration of the algorithm in \autocite{owen2011fast} shows that the same algorithm works in $\tau$-space and hence implies that geodesics are computable in polynomial time in $\tau$-space.
Indeed, geodesics in $\tau$-space satisfy the Characterisation Theorem \autocite[2.3--2.5 in][]{owen2011fast}, so once two $\tau$-trees are converted into an incompatibility graph, \autocite{owen2011fast} gives a polynomial algorithm to find the splits of the set of vertices of the graph corresponding to the geodesic.
However, the data structures for the algorithm in $\tau$-space should be different.
Indeed, the notion of compatibility is given in $\tau$-space by the notion of refinement (see Section~\ref{geoDataStrSection}), which encodes ranks in the incompatibility graph.
We have implemented the algorithm in Java and the implementation can be accessed at \autocite{tauGeodesic} under the GNU General Public Licence.
The novel data structures necessary for this implementation are formally introduced later in this paper.
The running time of the implementation is similar to that of \autocite{owen2011fast}.

Thus, we suggest that $\tau$-space serves as a tool for statistical analysis of stochastic processes over ultrametric phylogenetic trees.
Particularly, for computing the summary tree of a posterior sample obtained using, for example, MCMC.

\section{$\t$-space}\label{tSpace}

As absolute divergence times are often the object of interest, the parameterisation of trees using the times of their nodes is natural for several phylogenetic modes, e.g.\ birth-death models \autocite{kendall}.
Furthermore, birth-death priors are one of the main classes of priors used in Bayesian inference.
The purpose of this section is twofold.
First, we would like to study geometric and efficiency properties of one of the prominent parameterisations in evolutionary biology.
Second, we demonstrate how radically these properties can change after a seemingly negligible change in parameterisation.
Namely, converting $\tau$-coordinates to their initial sums, that is, to the absolute times of divergence events, makes fundamental results from combinatorial geometry such as Gromov's theorem used to prove Theorem~\ref{uniqueGeodesics} inapplicable, along with the algorithmic results from \autocite{owen2011fast} and \autocite{Owen2011-cc}.

We note that \textcite{Ardila2006-ov} considered branch lengths and a tree height to parameterise the space of ultrametric trees and prove that ultrametric trees are in one-to-one correspondence with the Bergman fans of complete graphs \autocite[for details see][]{Ardila2006-ov}.
As explained in Section~\ref{preAttemptSection}, this parameterisation is not convenient for our purposes.
However, the results in \autocite{Ardila2006-ov} imply that $\t$-space is in one-to-one correspondence with the Bergman fans of complete graphs, so our results on $\t$-space are applicable to the corresponding fans.

\subsection{Construction of space and uniqueness of geodesics}

Let us consider a completely resolved ultrametric tree $T$ with ranked topology $\rt(T)$ with no pair of nodes of the same rank.
For each node $v_i$ from $T$, let $t_i$ be the distance from $v_i$ to the nearest taxon.
In this way, we assign times to all nodes of $T$, with all taxa being of time $0$.
Let us order all internal vertices of $T$ according to their times: $v_1,\ldots,v_{n-1}$ (the ordering is the same as the ordering according to their ranks in $\rt(T)$).
Then the coordinates of tree $T$ in $\t$-space is the tuple $(\rt(T),t_1,\ldots,t_{n-1})$.

\begin{figure}
\centering
\includegraphics[width=0.6\textwidth]{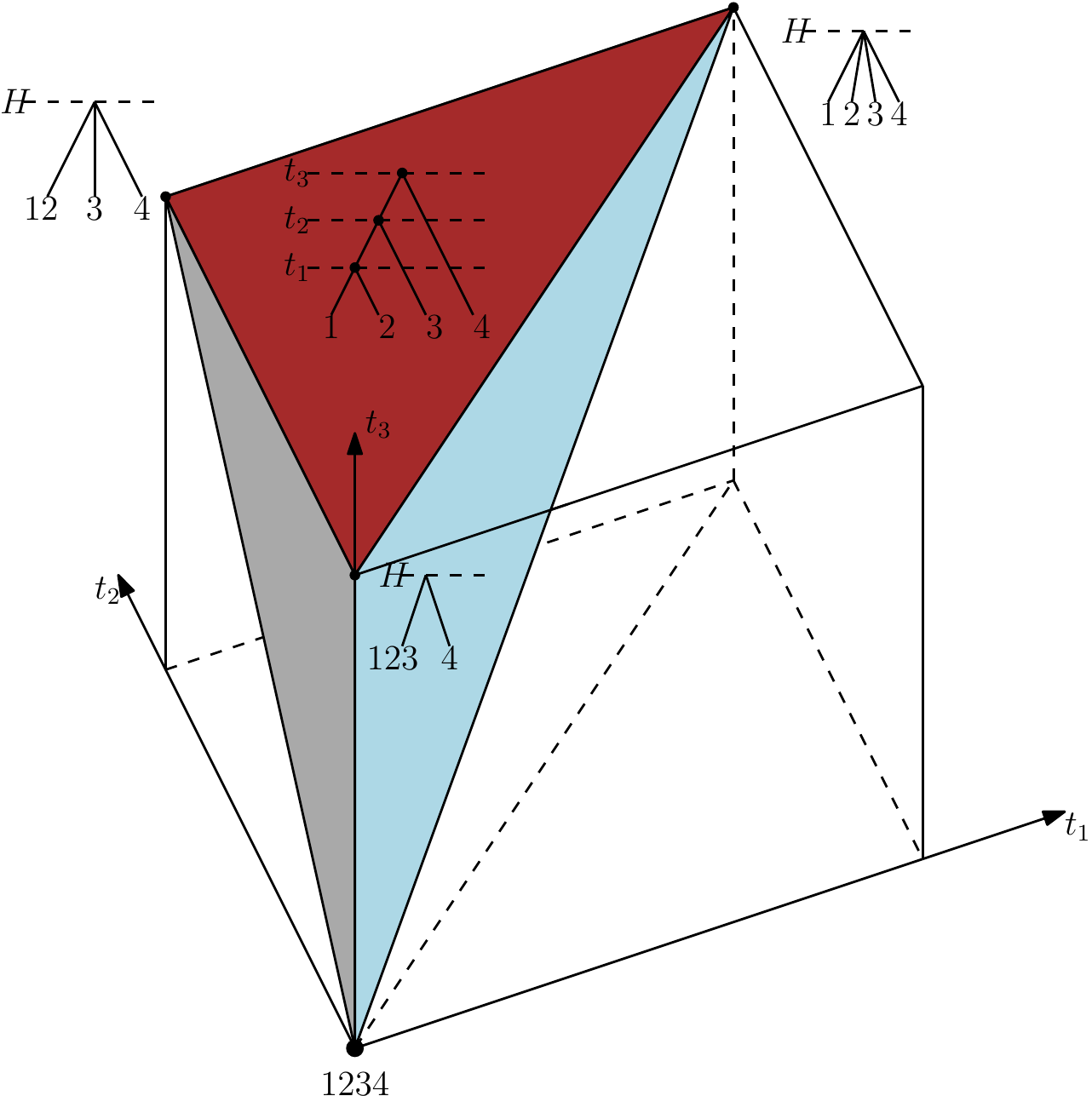}
\caption{One simplex of $\t$-space.}
\label{simplex.pdf}
\end{figure}

We note that if we vary the times of the nodes of $T$ while keeping the ranked topology preserved, we get a simplex $\{(t_1,\ldots,t_{n-1})\mid 0\leq t_1\leq\ldots\leq t_{n-1}\leq H\}$, where $H$ is an (artificial, sufficiently large) upper bound on the height of the
tree\footnote{Again, all in this paper remain true in the unbounded case with no $H$.}.
Figure~\ref{simplex.pdf} depicts one such simplex.

\begin{figure}
\centering
\includegraphics[width=0.7\textwidth]{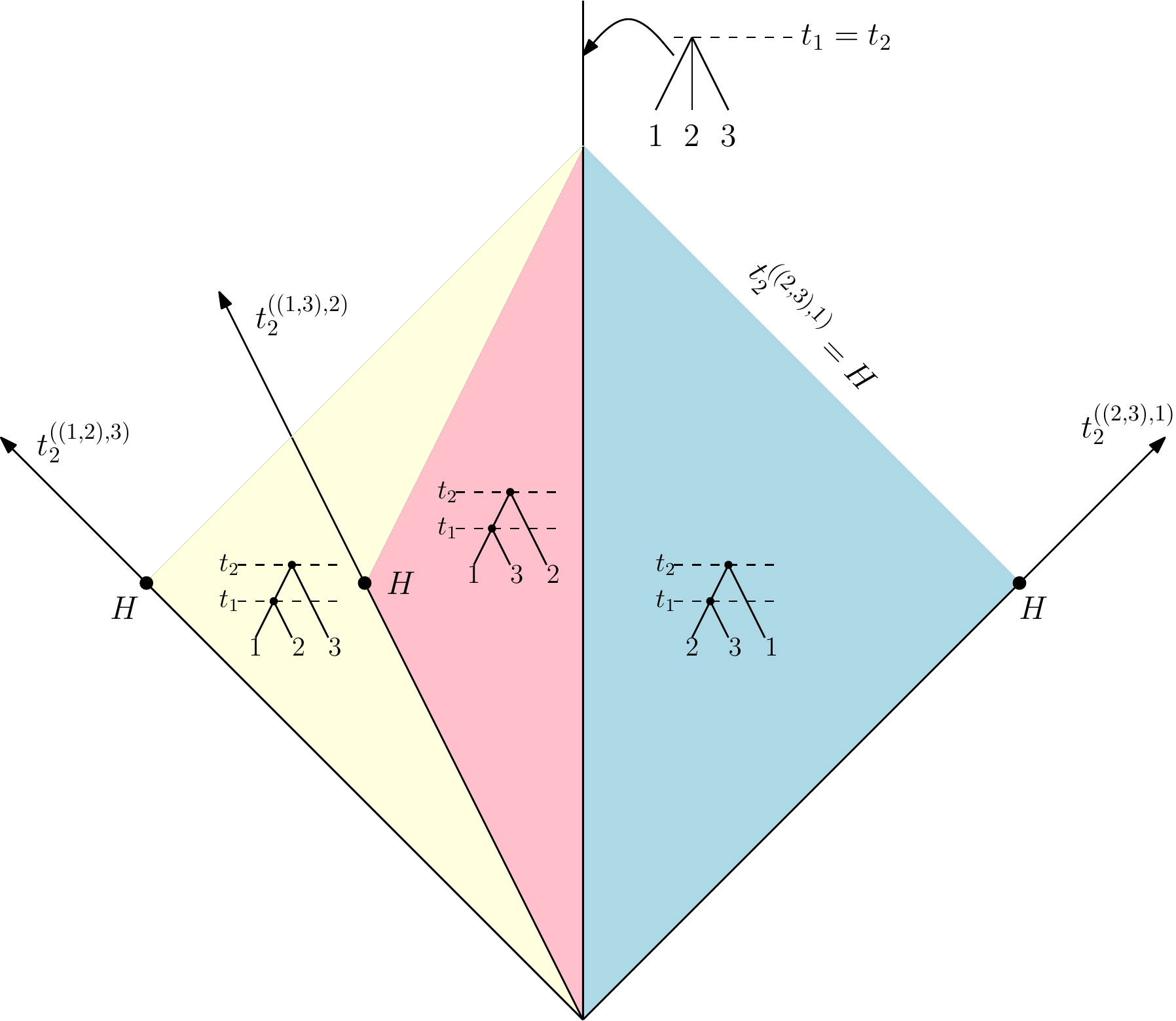}
\caption{Space $\TT_3$.
Three coloured triangles of the simplicial complex correspond to the three depicted topologies.
The triangles share a line that corresponds to the unresolved tree on three taxa.
The upper bound of the triangles is the artificial bound $H$.}
\label{tSpace2d.pdf}
\end{figure}

\begin{figure}
\centering
\includegraphics[width=0.7\textwidth]{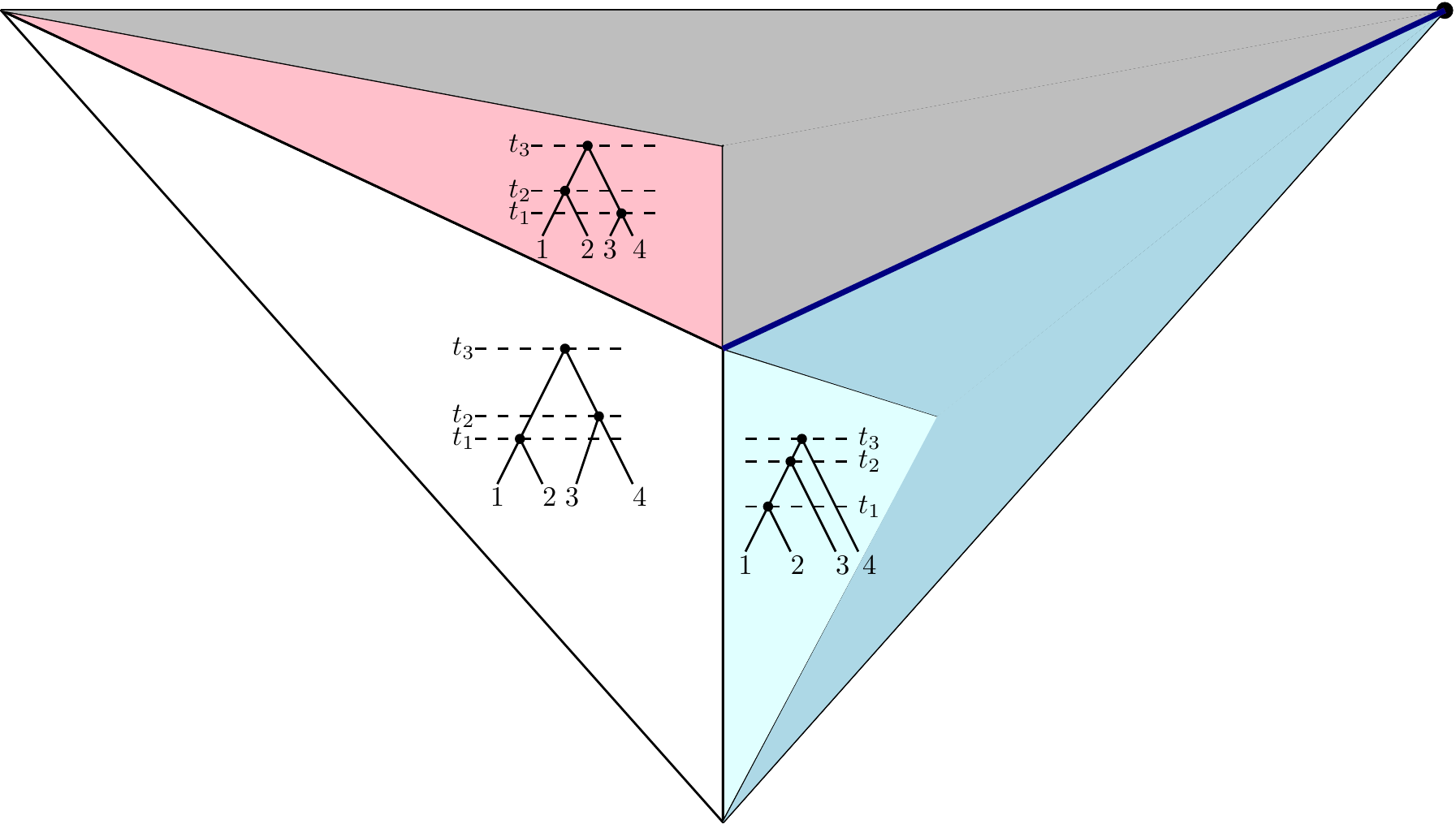}
\caption{One-sixth of 4D $\t$-space $\TT_4$ corresponding to the depicted topologies.
Unlike in Figure~\ref{tauSpace.pdf}, the simplices are not projected onto a $2$-dimensional subspace and drawn as $3$-dimensional pyramids.
Three such pyramids are depicted in white, blue, and grey.
The white pyramid shares a facet with both grey and blue pyramids.
The grey and blue pyramids share the edge $t_1=t_2=t_3$ of the complex only.}
\label{tSpace.pdf}
\end{figure}

We create a simplicial complex out of $\frac{(n-1)!\cdot n!}{2^{n-1}}$ such simplices corresponding to different ranked topologies on $n$ taxa in a similar way the complex is created in $\tau$-space, namely, we identify faces of simplices corresponding to the same tree topology.
The metric is defined in the same way as in $\tau$-space to be the standard piecewise Euclidean distance.
The first substantial difference is that the edge of the complex that is shared by all simplices is not an axes but rather the line $t_1=\ldots=t_{n-1}$.
Furthermore, the faces are defined by some of the coordinates being equal, $t_i=t_k$, rather than some of the coordinates being $0$.
We call the space so defined a {\em $\t$-space} and denote the $t$-space on $n$ taxa by $\TT_n$.
Figure~\ref{tSpace2d.pdf} depicts space $\TT_3$ in full and Figure~\ref{tSpace.pdf}---a part of $\TT_4$.

The following lemma is an important property that connects $\BHV$ space, $\tau$-space, and $\t$-space.
This lemma has (implicitly) been used to prove that $\tau$-space is $\CAT(0)$.

\begin{lemma}\label{homeomorphic}
$\BHV$ space, $\tau$-space, and $\t$-space are pairwise homeomorphic.
\end{lemma}

\proof
The homeomorphisms are induced by parameterisations $p_{\BHV}$, $p_{\tau}$, and $p_\t$ used to construct the corresponding spaces.
\endproof

In particular, this lemma implies that both $\tau$-space and $\t$-space are connected and simply connected.

Our next step is naturally to ask whether $\t$-space has unique geodesics.
This question cannot be answered in the same way as it is done for $\BHV$ and $\tau$-space using Gromov's characterisation of $\CAT(0)$ cubical complexes because $\t$-space is not a cubical complex.
Several sufficient conditions are known for simplicial complexes to be $\CAT(0)$ \autocite{januszkiewicz2006simplicial,fisher2007geometry}, however $\t$-space does not satisfy those conditions and hence they cannot be used to prove that $\t$-space is $\CAT(0)$.
Below we prove that $\t$-space satisfies the $\CAT(0)$ axiom thus giving a new and important example of a $\CAT(0)$ simplicial complex.

\begin{theorem}\label{tSpaceCAT0}
$\t$-space has unique geodesics.
\end{theorem}

\proof
We call a facet {\em shared} if the facet belongs to at least two different simplices in $\t$-space.

\begin{lemma}\label{tSpaceAngles}
Let $S$ be a simplex in $\TT_n$ and $\angle$ be an angle between a pair of shared facets in $S$.
Then $\angle \geq \pi/3$.
\end{lemma}

\proof
First, scale the simplex $S$ so that the height of the trees corresponding to $S$ is bounded by $1$.
Then the set of vertices $V$ of simplex $S$ is
\[
\{(\underbrace{0, \ldots, 0}_{n-1-i}, \underbrace{1, \ldots, 1}_i) \mid i\leq n-1\}.
\]
The set of facets of simplex $S$ is then given by the set of all $(n-1)$-element subsets of $V$.
Note that there are exactly two facets of $S$ which are not shared by any other simplex in $\TT_n$.
Indeed, the facet given by the set $V \setminus \{(0,\ldots,0)\}$ corresponds to a completely resolved topology and belongs to exactly one simplex---the one that corresponds to that topology.
The facet given by the set $V \setminus \{(1,\ldots,1)\}$ corresponds to a topology where a pair of taxa $a,b$ are unresolved at present and the rest of the nodes of the topology are resolved.
Since there is only one possibility to resolve the degenerate cherry $(a,b)$, this facet belongs to exactly one simplex as well.
Hence the shared facets of $S$ are precisely those that contain both $(0, \ldots, 0)$ and $(1, \ldots, 1)$.
For example, the grey and the red facets of the simplex in Figure~\ref{simplex.pdf} are not shared by any other simplex, while the blue and the invisible facets are shared with other simplices.

For every pair of shared facets of $S$, we now find the angle between them.
To do so we first need to find the normal vectors of all shared facets.
Every shared facet $F_i$, $1 \leq i \leq n-2$, of $S$ can be represented as an $(n-1) \times (n-1)$ matrix the $j$-th row of which is
\[
\begin{aligned}
(0, \ldots, 0) & \mbox{ ~~~ if ~~~ } j = 0,\\
(\underbrace{0, \ldots, 0}_{n-1-j}, \underbrace{1, \ldots, 1}_j) & \mbox{ ~~~ if ~~~ } 1 \leq j < i,\\
(\underbrace{0, \ldots, 0}_{n-1-(j+1)}, \underbrace{1, \ldots, 1}_{j+1}) & \mbox{ ~~~ if ~~~ } i \leq j < n-2,\\
(1, \ldots, 1) & \mbox{ ~~~ if ~~~ } j = n-2.
\end{aligned}
\]
In other words, facet $F_i$ is defined by removing the row $(\underbrace{0, \ldots, 0}_{n-1-i}, \underbrace{1, \ldots, 1}_i)$ from the set $V$.
Normal vectors $f_i$ of facets $F_i$ are then given by the null space of these matrices.
For every $i$ such that $1 \leq i \leq n-2$, we fix one such vector $f_i$ to be $(\underbrace{0, \ldots, 0}_{n-2-i}, (-1)^{s_i}, (-1)^{s_{i}+1}, \underbrace{0, \ldots, 0}_{i-1})$, where $s_i \in \{0,1\}$ is chosen depending on from what side of $F_i$ the simplex is located.

Since the inner product of pairs of these vectors is $-1$, $0$, or $1$, the smallest possible angle between the facets is $\pi/3$, which proves the lemma.
\endproof

We now note the following property.
For every positive real number $\varepsilon$, the $\varepsilon$-neighbourhood of the
origin\footnote{This is true not only for the origin but for every $\varepsilon$-neighbourhood of every star-tree.}
of $\TT_n$---point $(0,\ldots,0)$ in terms of the proof of previous lemma---contains a simplicial complex similar to $\TT_n$.
This property follows by scaling $\TT_n$ with small enough scaling factor.
Hence, to establish the claim of the theorem it is enough to show that the space is locally $\CAT(0)$.
To prove this, we apply the following characterisation of locally $\CAT(0)$ simplicial complexes:

\begin{theorem}[\cite{gromovOriginal}, {\cite[see][5.2 on p.\,206]{bridsonBook}}]
A finite simplicial complex is locally $\CAT(0)$ if and only if the link of every vertex of the complex is a $\CAT(1)$ space.
\end{theorem}

Hence, to finish the proof of Theorem~\ref{tSpaceCAT0}, we need to show that the link of every vertex of $\TT_n$ is a $\CAT(1)$ complex.
By Theorem 5.4(7) in \autocite[][p.\,206]{bridsonBook}, it is enough to show that $\TT_n$ contains no isometrically embedded circles of length less than $2\pi$.
That means the following.
Let $C$ be a geodesic curve in (the link of a vertex of) $\TT_n$ of length $\ell$ which is isometric to a Euclidean circle $C_E$.
If we scale the space so that $C_E$ is a unit circle then $\ell \geq 2\pi$.

This last property follows from Lemma~\ref{tSpaceAngles} for space $\TT_n$.
Let $T_1,\ldots,T_k$ be completely resolved pairwise different ranked tree topologies with the property that if $C$ intersects a simplex corresponding to tree topology $T$ then there is an $i$ such that $T = T_i$.
Furthermore, we assume that the topologies $T_1,\ldots, T_k$ are ordered as they are traversed by $C$.
Since the length $\ell$ of $C$ satisfies $\ell \geq k*\angle$ where $\angle$ is the smallest angle between facets of the simplex corresponding to a $T_i$, the cases when $k \geq 6$ follow from Lemma~\ref{tSpaceAngles} directly.

To finish the proof, we consider the cases when $k \leq 5$.
Since the shortest possible cycle is of length $4$, it remains to consider only two cases:
\begin{itemize}
\item[$k = 4$.] Let $F_{i_1}, F_{i_2}, F_{i_3}, F_{i_4}$ be the facets represented in the form of
matrices\footnote{It is important to note that $i = j$ does not imply $F_i = F_j$, as these could correspond to different tree topologies.
However, this notation is convenient and does not result in ambiguities in this proof.}
as above in the order they are crossed by the geodesic circle $C$.
That is, $F_{i_1}$ is the facet shared by $T_1$ and $T_2$, $F_{i_2}$---by $T_2$ and $T_3$, $F_{i_3}$---by $T_3$ and $T_4$, and $F_{i_4}$---by $T_4$ and $T_1$.
Recall that indices $i_1$, $i_2$, $i_3$, and $i_4$ correspond to $\t$-coordinates on which the corresponding topology (rank) move is performed.

First, assume that $|i_1 - i_2| > 1$.
In this case, $i_3 = i_1$ and $i_4 = i_2$, or $i_3 = i_2$ and $i_4 = i_1$.
The latter case is not possible as in that case $F_{i_3} = F_{i_2}$ and hence $T_2 = T_4$.
Since $|i_1 - i_2| > 1$, the scalar product of normal vectors corresponding to $F_{i_1}$ and $F_{i_2}$ is $0$ and the angle between them is $\pi/2$.
Hence, all four angles are $\pi/2$ each and hence $\ell = 2\pi$.

Now assume $|i_1 - i_2| = 1$.
In this case the cycle cannot exist.
Indeed, one has to consider all possible combinations of tree topology changes (moves) corresponding to $F_{i_1}$ and $F_{i_2}$: both are rank changes (RR), the first is an $\NNI$ move and the second is a rank change (NR), similarly RN and NN.
In all the four cases the circle has to cross both $F_{i_1}$ and $F_{i_2}$ twice and hence $T_2$ has to coincide with $T_4$.
We consider the NN case and the other cases follow similarly.
If both $F_{i_1}$ and $F_{i_2}$ correspond to different $\NNI$ moves (as in NN), then tree $T_3$ is at $\NNI$ distance $2$ from $T_1$.
In this case, the only $\NNI$ path from $T_3$ to $T_1$ of length $2$ has to follow the moves corresponding to $F_{i_2}$ and $F_{i_1}$.

\item[$k = 5$.] Let $F_{i_1}, F_{i_2}, F_{i_3}, F_{i_4}, F_{i_5}$ be facets as above.
We show that the cycle cannot exist.
Let us consider the types of tree topology moves corresponding to $F_{i_1}$, $F_{i_2}$, and $F_{i_3}$.
First, consider the case when all three moves are rank moves, RRR.
In this case, three nodes have changed their ranks and at least three rank changes are necessary to return to the original tree $T_1$, hence the cycle of length $5$ is not possible.
Second, consider the case when exactly one of the tree moves is an $\NNI$ move.
In this case, the tree obtained after the three steps is at $\NNI$ distance $1$ plus two nodes have changes their ranks. Again, at least three steps are necessary to return to tree $T_1$.
Third, consider the case when exactly two of the tree moves are $\NNI$.
In this case, the tree obtained is at $\NNI$ distance two from $T_1$ and one node changes its rank.
Again, at least two $\NNI$'s and one extra move are necessary to return to $T_1$.

Finally, let us assume that all facets $F_{i_1}$, $F_{i_2}$, and $F_{i_3}$ correspond to $\NNI$ moves, NNN.
Note that if after the three $\NNI$ moves we obtain a tree at $\NNI$ distance $3$ from the original tree $T_1$, then the cycle of length $5$ cannot exist.
Hence we assume that the $\NNI$ distance between $T_1$ and $T_4$ is $2$.
In this case, all the five facets correspond to $\NNI$ moves, and that is not possible as at least one rank move is necessary.
Indeed, after the $\NNI$ move corresponding to $F_{i_1}$, there must be a node that changes its rank in the following sense: there exist taxa $A$ and $B$ such that the difference of ranks of their most recent common ancestors $\mrca(A, B)$ in $T_1$ and $T_2$ is $1$.
In order for the cycle to return to $T_1$ without intersecting $F_{i_1}$ twice, there must exist a rank move on the way back to $T_1$.
\endproof
\end{itemize}

The change in the parameterisation of trees results not only in the question of uniqueness becoming more complicated.
What is also important is that the algorithms used for computing geodesics in $\BHV$ and $\tau$-space cannot directly be applied in $\t$-space.
Moreover, their existence has to be questioned.
Hence we propose the following problem, on which we make some progress below but do not obtain a complete answer.

\begin{problem}\label{tSpaceAlgorithm}
What is the complexity of computing geodesics in $\t$-space?
\end{problem}

\subsection{Geometry and data structures}\label{geoDataStrSection}

\begin{figure}
\centering
\includegraphics[width=0.7\textwidth]{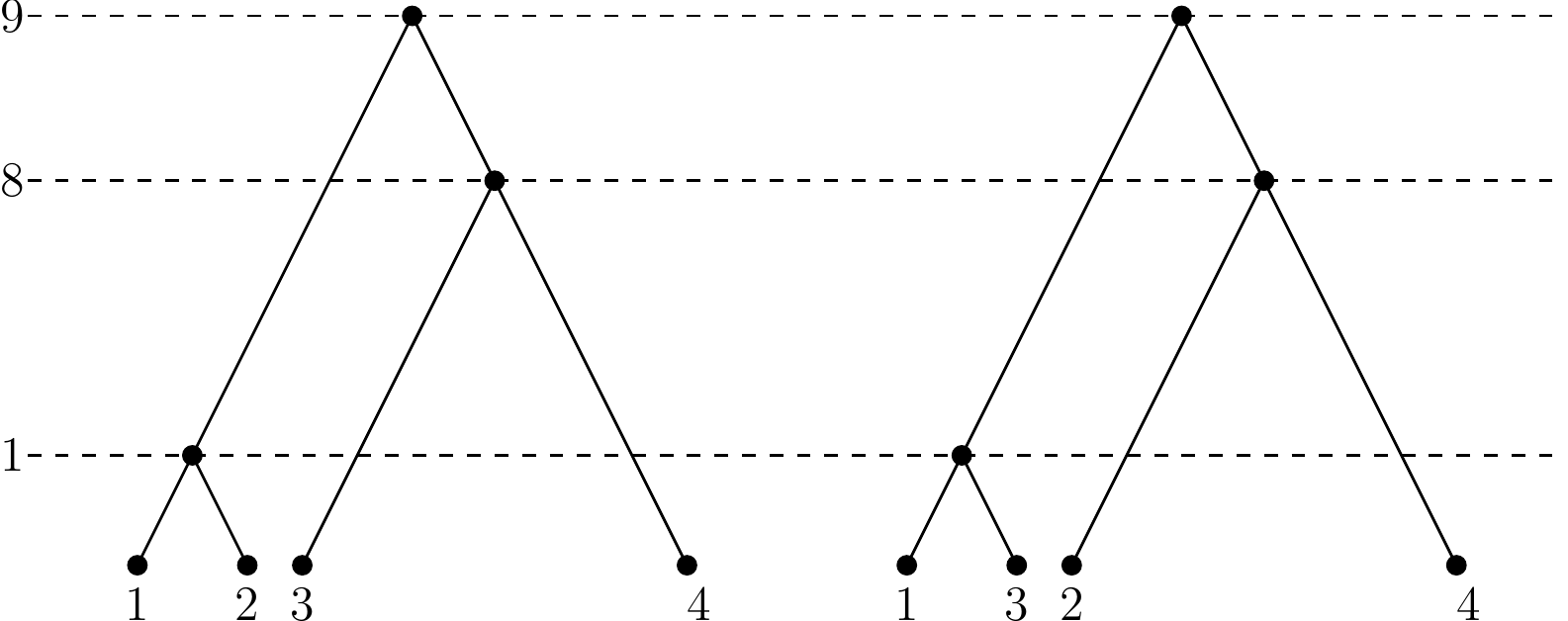}
\caption{$\BHV$- and $\tau$-geodesics are cone-paths while $\t$-geodesic is not.}
\label{conePaths.pdf}
\end{figure}

One of the key properties that make $\tau$-space and $\t$-space so different is that the cone-path is rarely a geodesic in $\t$-space.
Indeed, in both $\BHV$ and $\tau$-space the position of two cubes can result in the cone-path being the geodesic between every pair of trees from those cubes.
For example, if two trees $T$ and $R$ have topologies with no compatible splits then the geodesic between $T$ and $R$ is a cone-path \autocite{bhv}.
$\t$-space does not have this property.
Let us illustrate this effect with the following example.
Consider trees $T$ and $R$ depicted in Figure~\ref{conePaths.pdf}.
Since the trees do not have compatible splits, the geodesic is a cone-path in both $\BHV$ and $\tau$-space.
However, the shortest cone-path in $\t$-space passes through the star-tree of height $6$ and has length $2\sqrt{38}>12.3$, while the path that goes through $((1,2)\colon 4,3,4)\colon 8$, $((1,3,2)\colon 6,4)\colon 8$, and $((1,3)\colon 4,2,4)\colon 7$ (the numbers following the colon are heights of the corresponding clades) has length $\sqrt{10}+\sqrt{8}+\sqrt{6}+\sqrt{14}<12.2$ and is hence shorter than every cone
path\footnote{However, this path is not a geodesic either.
To find the actual geodesic between these two trees is a simple but interesting exercise.}.

This example demonstrates another important property that distinguishes $\t$-space from $\tau$-space.
Every tree on the geodesic between two trees in $\tau$-space contains only splits that present in the origin tree or in the destination tree (or both).
Not so in $\t$-space.
Split $(123 \mid 4)$ does not present in either tree in Figure~\ref{conePaths.pdf} but present in an intermediate tree on the $\t$-geodesic.

We now develop a formalism that is convenient for the study of the computational content of $\t$-space.
This formalism can be used for $\tau$-space as well and is actually the data structure that we use in our implementation of geodesic algorithms for $\tau$-space \autocite{tauGeodesic}.
The formalism is motived by and consistent with the treatment of ranked tree topologies in \autocite{steel}.

By a {\em partition with attached time coordinate}, we mean an object of the form $(N_1\mid\ldots\mid N_q)\colon t$, where $N_1\mid\ldots\mid N_q$ is a partition of taxa that can be obtained by cutting the tree along the line obtained by fixing the time coordinate, and $t$ is the least value of the time coordinate that produces this partition.
For example, the left-hand side tree in Figure~\ref{conePaths.pdf} is defined by the set of partitions
\[
\{(12\mid3\mid4)\colon 1,~(12\mid34)\colon 8,~(1234)\colon 9\},
\]
while the right-hand side tree---by
\[
\{(13\mid2\mid4)\colon 1,~(13\mid24)\colon 8,~(1234)\colon 9\}.
\]
Removing one or more partitions from a set of partitions that defines a completely resolved tree gives rise to a non-resolved tree or a tree with two or more internal nodes of the same rank.
For example, if we remove partition $(12\mid3\mid4)$ from the left-hand side tree in Figure~\ref{conePaths.pdf} then we get the tree $((1,2),(3,4))$ of height $9$ with both internal nodes being of height $8$.
Alternatively, if we remove partition $(12\mid34)$ from that tree then we get the unresolved tree $((1,2),3,4)$ of height $9$ with a common ancestor of $(1,2)$ at height $1$.

We note here that as we consider only trees with all taxa being at time $0$, the partition $(1\mid 2\mid\ldots\mid n)\colon 0$ is assumed to be (invisibly) present
everywhere\footnote{This assumption cannot be made in the general setting of time-trees or even more general setting of sampled ancestor trees \autocite{gavryushkina2013recursive,gavryushkina2014bayesian}.}.
Clearly, a tree is unambiguously defined by its set of partitions with attached time coordinates, and a set of partitions defines a tree if and only if one member of every pair of partitions from the set {\em refines} the other and the time coordinates of the partitions are monotonic under those refinements.

The fact that the restriction of a geodesic to a simplex is a straight line justifies the following definition.

We assume that trees $T$ and $R$ are completely resolved and have all internal nodes of different ranks, that is, neither of them has $\t$-coordinates $t_i$ and $t_{i+1}$ such that $t_i = t_{i+1}$.
We say that the geodesic $\gamma$ between trees $T$ and $R$ is {\em computable (in polynomial time)} if (a polynomial and) an algorithm exists that given the sets of partitions with attached time coordinates $A_T$ and $A_R$ that define $T$ and $R$ respectively, outputs (after a number of steps bounded by the polynomial of $n$) a sequence of sets of partitions $A_0, \ldots, A_k$ with time coordinates attached to every partition such that the following two properties are satisfied:

\begin{itemize}
\item $A_T = A_0$ and $A_R = A_k$.
\item For every $i < k$, the pair of sets $A_i, A_{i+1}$ along with the attached time coordinates defines the trees where geodesic $\gamma$ enters and exits simplex $S_i$, respectively.
Here, $S_0, \ldots, S_{k-1}$ are all the simplices geodesic $\gamma$ traverses in the order they are traversed.
Particularly, $S_0$ contains $T$ and $S_{k-1}$ contains $R$.
\end{itemize}

Since tree $T$ ($R$) is completely resolved, the number of elements of $A_0$ ($A_k$) is $n-1$.
In terms of simplices, the number $n - 1 - |A_i|$, where $|A_i|$ is the number of elements in $A_i$, is the codimension of the face of simplex $S_i$ where $\gamma$ enters $S_i$.
In terms of trees, this number is the number of multifurcations plus the number of non-resolved ranks of internal nodes of the tree corresponding to $S_i$.

Note that the properties above imply that all sets $A_i$ are pairwise different, all time coordinates attached to the partitions from the same set are pairwise different, and the time coordinates attached to the same partition in different sets may or may not be different.
Clearly, every geodesic is unambiguously defined by a sequence of sets of partitions with attached time coordinates satisfying these properties.

We now give an example of calculating a geodesic in $\t$-space.

\begin{example}
What is the geodesic between trees $E$ and $S$ depicted in Figure~\ref{NNIpathLemma.eps}?

\begin{figure}
\centering
\includegraphics[width=.63\textwidth]{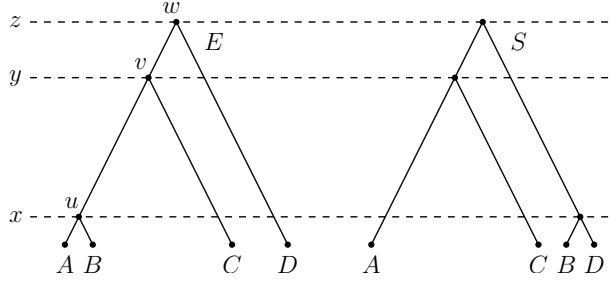}
\caption{The geodesic between trees $E$ and $S$ is or is not a cone-path depending on the height of the root.
$A$, $B$, $C$, and $D$ are some trees and not necessarily taxa.}
\label{NNIpathLemma.eps}
\end{figure}

Direct computations show that the geodesic between the two trees is a cone-path if and only if $y \geq \frac{x + z}{2}$, in which case the geodesic passes through the star-tree of height $\frac{x+y+z}{3}$.
If $y < \frac{x + z}{2}$ then the geodesic passes through the tree $\{(AC\mid B\mid D)\colon y, (ABCD)\colon \frac{x + z}{2}\}$, that is, node $v$ does not move along the geodesic.
Since geodesics restricted to a simplex are straight lines, this provides a complete characterisation of $\t$-geodesics between pairs of trees having the topologies of $E$ and $S$.

The (inefficient) brute-force algorithm for computing these geodesics would be to consider (not necessarily shortest) paths in the $\NNI$ graph that lead from one tree to the other and solve the minimisation problem to obtain the $\t$-coordinates of the nodes of unresolved trees, that is, the coordinates of the intersection points on faces.
\end{example}

This example is an illustration of the following fundamental property of $\t$-space.
Recall that in both $\BHV$ and $\tau$-space there exist pairs of tree topologies such that the geodesic between trees having those topologies is a cone-path no matter what the branch lengths are.
We now show that $\t$-space does not have this property.

\begin{theorem}
For every pair of ranked tree topologies $\rt_1$ and $\rt_2$, there exist trees $T$ and $R$ such that $\rt(T) = \rt_1$, $\rt(R) = \rt_2$, and the geodesic $\gamma$ between $T$ and $R$ is not a cone-path.
\end{theorem}

\proof
Assume first that for every pair of taxa $s,k$ at least one of the nodes $\mrca_{\rt_1}(s,k)$, $\mrca_{\rt_2}(s,k)$ is a root node.
This can only happen when the number of taxa is three or the number of taxa is four and the topologies of the trees are of the form $((1,2),(3,4))$ and $((1,3),(2,4))$.
In the first of these cases, the theorem is trivially true.
Consider the second case.
Between every pair of completely resolved trees with topologies $((1,2),(3,4))$ and $((1,3),(2,4))$ there exists a path that crosses exactly $3$ facets.
Since the angle between those facets is $\pi/3$, the trees $T$ and $R$ can be chosen so that the angle between them is less than $\pi$.
Hence thus chosen trees $T$ and $R$ can be connected by a straight line, which is the shortest path possible.

Assume now that there exists a pair of taxa $s,k$ such that both $\mrca_T(s,k)$ and $\mrca_R(s,k)$ are not root nodes.
Let $T$ and $R$ be an arbitrary pair of trees such that $\rt(T) = \rt_1$ and $\rt(R) = \rt_2$.
We prove the following stronger version of the theorem:

{\em Let $T_\delta$ be the tree obtained from $T$ by increasing the root height by $\delta$.
Then there exists a number $H$ such that a path where $\mrca_S(s,k)$ is not a root node for all trees $S$ on the path is shorter than every cone-path from $T_\delta$ to $R_\sigma$, where $\delta, \sigma \geq H$.
Hence the geodesic between $T_\delta$ and $R_\sigma$ cannot be a cone-path.
Furthermore, the number $H$ is computable in polynomial time from $T$ and $R$.}

Let $\{t_i\}$ and $\{r_i\}$ be the time coordinates of $T$ and $R$ respectively.
Then the shortest cone-path from $T$ to $R$ passes through the star-tree of height $h$ obtained from the following minimisation:
\begin{equation}
\sqrt{\sum\limits_{i=1}^{n-1}(h - t_i)^2} + \sqrt{\sum\limits_{j=1}^{n-1}(h - r_j)^2} \to \min
\tag{$*$}
\label{coneThEqn}
\end{equation}

Hence, the height of the star-tree can be made arbitrarily high by increasing the heights $t_{n-1}$ and $r_{n-1}$ of trees $T$ and $R$.

Fix a large enough number $H$ (the exact value will be determined later) obtained from minimisation~(\ref{coneThEqn}) and consider the following path between $T_\delta$ and $R_\sigma$, where $\delta,\sigma \geq H$.
Fist, the path follows the straight line from $T_\delta$ to the tree $S'$ that has two internal nodes: $\mrca_{S'}(s,k)$ and the root, the time coordinate of $\mrca_{S'}(s,k)$ equals to the time coordinate of $\mrca_T(s,k)$ and the time coordinate of the root equals to $H$.
Then, the path follows the straight line from $S'$ to $S''$, where the only difference between $S'$ and $S''$ is that the time coordinate of $\mrca_{S''}(s,k)$ equals to the time coordinate of $\mrca_R(s,k)$.
Finally, the path follows the straight line from $S''$ to $R$.
Let $i_s$ be the number of the time coordinate of $\mrca_T(s,k)$ and $j_k$---of $\mrca_R(s,k)$
Then the length of this path is equal to
\begin{equation}
\sqrt{\sum\limits_{i\ne i_s}(H - t_i)^2} + \sqrt{\sum\limits_{j\ne j_k}(H - r_j)^2} + |r_{j_k} - t_{i_s}|
\tag{$**$}
\label{coneThEqn2}
\end{equation}

Note that the value of this function~(\ref{coneThEqn2}) is smaller than the value of the objective function in minimisation~(\ref{coneThEqn}) for all large enough values of $H$.
This is our first requirement on the number $H$.
The second requirement is that the path described in~(\ref{coneThEqn2}) exists.
Since both these requirements can be checked in polynomial time, the stronger version of the theorem is proved.
\endproof

It follows from the proof of this theorem that for every pair of simplices in $\t$-space, a non-trivial part of them consists of trees between which the geodesic is not a cone-path.
Hence, the volume of pairs of trees between which the geodesic is not a cone-path is positive for every pair of simplices, unlike in $\BHV$ or $\tau$-space.
This is because geodesics in $\t$-space often follow $\NNI$-paths.
However, they do not necessarily follow {\em shortest} $\NNI$-paths.
Consider two caterpillar trees $(((((((((1,2),3),4),5),6),7),8),9),10)$ and $(((((((((1,2),5),6),7),8),9),3),4),10)$.
It is not hard to see that no tree, apart from these two, is a caterpillar tree on an $\NNI$-geodesic between them, however every tree on the geodesic between trees with these tree topologies is a caterpillar tree in $\t$-space.
A more detailed investigation of this phenomenon is the subject of our future work.

\section{Conclusion and further directions}

We have considered two standard parameterisations of the space of ultrametric phylogenetic trees: (1) using lengths of coalescent intervals and (2) using times of divergence events.
By considering suitable polyhedral complexes, we have found two possible representations of the space of trees called $\tau$-space and $\t$-space respectively.
Despite their similarity, the two parameterisations have significantly different geometric and algorithmic properties.
For example, we showed that geodesics, and hence Fr\'echet means, are different in the two spaces.
Although it required quite different geometric approaches, we proved that shortest paths are unique in both $\tau$- and $\t$-space.
We also proved that shortest paths are efficiently computable in $\tau$-space.
We have implemented the algorithm for computing exact shortest paths in Java.
We also implemented the algorithms for efficiently approximating Fr\'echet means, standard variances, and some other geometric and statistical characteristics of finite samples of trees.
Although the algorithmic complexity of $\t$-space remains unknown, the space has several properties that are desirable for statistical analysis of tree space.
For instance, we proved that the paths that traverse a star-tree are often shortest in $\tau$-space and are rarely shortest in $\t$-space.
This feature of $\t$-space is a desirable property for phylogenetic applications, and particularly for summarising posterior samples by a point estimate.
Indeed, one of the unpleasant features of $\BHV$ space is that parts of the summary tree are often the star-tree, when incompatible subtrees are supported by the posterior \autocite{Hotz2013-bz}.
This feature is a consequence of a fundamental geometric property of the space---for some pairs of trees the shortest path traverses a star-tree no matter what the branch lengths of the trees are.
Both $\BHV$ and $\tau$-space suffer from this feature, while $\t$-space is free from it.
Thus we expect summary trees produced using $\t$-space to be more informative and realistic.

Although all results about $\t$-space in this paper are presented for ultrametric trees, they can be extended to the set of all time-trees, as well as to the set of all sampled ancestor trees.
In the light of the work of \textcite{Sturm2002-da} on statistics over $\CAT(0)$ spaces, this makes $\t$-space a very promising candidate for the role of {\em the} parameterisation of phylogenetic time-trees and sampled ancestor trees.
However, the details of the extension as well as the question of efficiency remain for the future work.

An obvious direction of further research is to test our algorithms on simulated and real data sets, compare them with known algorithms, and suggest what extra formal properties of a parameterisation of the tree space are desirable.
As is suggested in our work, there are other possible ways that ultrametric tree space can be parameterised.
We have considered two obvious parameterisations and established that they are already quite different.
One can certainly come up with many other ways to parameterise the space.
The question arises:

\begin{problem}
Is there, in some sense, a single {\em optimal} parameterisation of the tree space?
If not, what is the class of acceptable parameterisations?
\end{problem}

Our paper suggests a number of directions for further theoretical investigations.
An important statistical question is

\begin{problem}
What parameterisation should be used for coalescent models?
Birth-death models?
Must the parameterisations used for these two types of models be different?
\end{problem}

This problem is especially intriguing in the light of work of \textcite{Stadler2015-jl}.
The first step towards an answer for this question is obviously to consider the coalescent and the birth-death priors in $\tau$- and $\t$-spaces.
Are these priors continuous in these spaces?
Can the distance between two trees be made a (simple) function of their prior probabilities?

Although much work has been done to investigate $\CAT(0)$ simplicial complexes, no satisfactory characterisation of the complexes is known \autocite{fisher2007geometry}.
Further research is needed with an eye towards effectiveness properties.
The problem in general is expected to be hard because even constructing non-trivial examples of $\CAT(0)$ simplicial complexes requires significant effort and only a few such examples are known \autocite{fisher2007geometry}.
In this paper, we have provided such an example---the $\t$-space.
Hence the following question, which we ask for $\t$-space, is also important for $\CAT(0)$ simplicial complexes in general.

\begin{problem}
Is there an efficient (in any sense) algorithm for computing shortest paths between trees in $\t$-space?
\end{problem}

As we have established in this paper, the measure (volume) of the set of pairs of trees between which the shortest path traverses a star-tree is positive in $\tau$-space and $\t$-space.
This measure is positive in $\BHV$ space \autocite{bhv} as well.
Hence the obvious question to understand the geometry of the space is to find this measure.
More precisely:

\begin{problem}
Let $\mu_n$ be the uniform measure on the set of pairs of trees on $n$ taxa between which the geodesic is a cone
path\footnote{We assume that the measure is the uniform probability measure, that is, the measure of the whole space is equal to $1$ and all cubes (simplices) have equal measures.}.
What is the value of $\mu_n$ for $\BHV$ space?
For $\tau$-space?
For $\t$-space?
Is the sequence $\{\mu_n\}_{n\in\omega}$ convergent?
If so, what is the limit $\lim_n\mu_n$?
What is the meaning of this limit?
\end{problem}

Clearly, $\mu_3=1$ in all $\BHV$, $\tau$- and $\t$-space.
To find $\mu_4$ is an entertaining exercise.

\section*{Acknowledgements}

We are grateful to Mike Steel, David Bryant, and Steffen Klaere for useful comments and advice on possible approaches to the problems considered in this paper.
We thank Alexandra Gavryushkina, Tim Vaughan, and St\'ephane Guindon for valuable advice on presentation of the results, and Dmitry Berdinsky and Jacek {\'S}wi{\k{a}}tkowski for navigating us through modern combinatorial geometry.
Finally, we thank the editor and two anonymous referees, who reviewed this work for the Journal of Theoretical Biology, for their very valuable comments on the earlier versions of the manuscript.

AJD was partially supported by a Rutherford Discovery Fellowship and Marsden Fund from the Royal Society of New Zealand.
AG was supported by the Rutherford Discovery Fellowship awarded to AJD and by Marsden Fund grant from RSNZ.

\printbibliography

\end{document}